\documentclass{article}
\usepackage{fullpage}
\usepackage{amsmath}
\usepackage{amssymb}
\usepackage{hyperref}

\usepackage{epsfig}
\usepackage{color}

\usepackage{amsthm}

\title{Decay of the Maxwell field on the Schwarzschild manifold}
\author{P. Blue \\
Max Planck Institute for Gravitational Physics \\
(Albert Einstein Institute)\\
Am M\"uhlenberg 1\\
D-14476 Golm\\
Germany\\
pblue@aei.mpg.de}  

\newtheorem{theorem}{Theorem}
\newtheorem{lemma}[theorem]{Lemma}

\newtheorem{remark}[theorem]{Remark}

\newcommand{\Reals}{\mathbb R}
\newcommand{\Naturals}{\mathbb N}

\newcommand{\GenericM}{{\mathcal{M}}}
\newcommand{\rs}{{r_*}}
\newcommand{\horifac}{(1-2M/r)}
\newcommand{\scri}{\mathfrak{I}}
\newcommand{\dt}{\partial_t}
\newcommand{\dr}{\partial_{\rs}}
\newcommand{\drOrig}{\partial_r}
\newcommand{\dtheta}{\partial_{\theta}}
\newcommand{\dphi}{\partial_{\phi}}
\newcommand{\dtt}{\partial_t^2}
\newcommand{\drr}{\partial_{\rs}^2}
\newcommand{\drrr}{\partial_{\rs}^3}
\newcommand{\slap}{\Delta_{S^2}}
\newcommand{\dAng}{\not\!\nabla}

\newcommand{\vecX}{X}
\newcommand{\vecY}{Y}
\newcommand{\vecZ}{Z}

\newcommand{\vecT}{T}
\newcommand{\vecR}{R}
\newcommand{\vecH}{\Theta}
\newcommand{\vecPh}{\Phi}
\newcommand{\vecTunit}{{\hat{T}}}
\newcommand{\vecRunit}{{\hat{R}}}
\newcommand{\vecHunit}{{\hat{\Theta}}}
\newcommand{\vecPhunit}{{\hat{\Phi}}}
\newcommand{\vecConf}{K}
\newcommand{\formConfg}{\underline{K}}

\newcommand{\formXg}{\underline{\vecX}}

\newcommand{\formZg}{\underline{\vecZ}}
\newcommand{\formTg}{\underline{\vecT}}
\newcommand{\formRg}{\underline{\vecR}}

\newcommand{\formTunitg}{\underline{\vecTunit}}
\newcommand{\formRunitg}{\underline{\vecRunit}}

\newcommand{\veclcoord}{{L}}
\newcommand{\vecncoord}{{N}}
\newcommand{\vecmcoord}{{M}}
\newcommand{\vecmbcoord}{{\bar{M}}}
\newcommand{\vecl}{l_\text{ex}}
\newcommand{\vecn}{n_\text{ex}}
\newcommand{\vecm}{m}
\newcommand{\vecmb}{\bar{m}}
\newcommand{\veclunit}{\hat{l}}
\newcommand{\vecnunit}{\hat{n}}
\newcommand{\veclP}{{\tilde{l}}}
\newcommand{\vecnP}{{\tilde{n}}}
\newcommand{\vecea}{e_{\isa}}
\newcommand{\veceb}{e_{\isb}}

\newcommand{\ia}{\alpha}
\newcommand{\ib}{\beta}
\newcommand{\ic}{\gamma}
\newcommand{\id}{\delta}
\newcommand{\ie}{\epsilon}

\newcommand{\isa}{A}
\newcommand{\isb}{B}

\newcommand{\CoDeriv}{\nabla}
\newcommand{\Lie}{\mathcal{L}}
\newcommand{\Extd}{\text{d}}
\newcommand{\forms}{\Omega}
\newcommand{\sfc}{\mathcal{S}}
\newcommand{\tensorA}{A}
\newcommand{\g}{\text{\bf{g}}}
\newcommand{\MaxF}{F}
\newcommand{\StressEnergy}{\text{\bf T}}
\newcommand{\MaxA}{A}
\newcommand{\Lagrangian}{L}
\newcommand{\GenMomentum}[1]{{}^{(#1)}P}
\newcommand{\GenMomentumFormg}[1]{{}^{#1}\underline{P}}
\newcommand{\Deformation}[1]{{}^{(#1)}\pi}
\newcommand{\GenEnergy}[1]{E_{#1}}
\newcommand{\insertion}[1]{\text{i}_{#1}}
\newcommand{\GeneralField}{\phi}
\newcommand{\GenManifold}{M}
\newcommand{\Action}{S}
\newcommand{\VarDeriv}{\delta}

\newcommand{\Islab}[2]{{[#1,#2]\times\Reals\times S^2}}
\newcommand{\Islice}[1]{{\{#1\}\times\Reals\times S^2}}
\newcommand{\intTslabInfinity}[1]{\int_{[#1,\infty)\times\Reals\times S^2}}
\newcommand{\dThreeVec}{d\nu}

\newcommand{\dThreenTunit}{\horifac^\frac12 r^2d\rs d^2\omega}
\newcommand{\dThreeHF}{\horifac r^2d\rs d^2\omega}
\newcommand{\dThreePOneD}{d\rs d^2\omega}
\newcommand{\dFourPOneD}{d\rs d^2\omega dt}
\newcommand{\dFourGen}{d^4x}
\newcommand{\dFourEx}{\horifac r^2 d\rs d^2\omega dt}

\newcommand{\intTslice}[1]{\int_{\Islice{#1}}}
\newcommand{\intTslab}[2]{\int_{\Islab{#1}{#2}}}
\newcommand{\intTBox}[3]{\int_{\{#1\}\times(#2,#3)\times S^2}}
\newcommand{\intTLCslice}[1]{\int_{\{#1\}\times[-(3/4)#1,(3/4)#1]\times S^2}}
\newcommand{\intTLCtighterslice}[1]{\int_{\{#1\}\times[-(1/2)#1,(1/2)#1]\times S^2}}
\newcommand{\intTLCintermediateslice}[1]{\int_{\{#1\}\times\{ |\rs|\in[(1/2)#1,(3/4)#1]\}\times S^2}}

\newcommand{\chiTrap}{\chi_{\text{trap}}}

\newcommand{\Elec}{\vec{E}}
\newcommand{\Mag}{\vec{B}}
\newcommand{\MyFp}{\phi_1}
\newcommand{\MyF}{\phi_0}
\newcommand{\MyFm}{\phi_{-1}}
\newcommand{\Myrho}{\rho}
\newcommand{\Mysigma}{\sigma}
\newcommand{\MyFi}{\phi_i}
\newcommand{\CKrho}{{\rho}}
\newcommand{\CKsigma}{{\sigma}}
\newcommand{\CKalpha}{{\alpha}}
\newcommand{\CKalphab}{{\underline{\alpha}}}
\newcommand{\CoordWtp}{\Phi_{1}}
\newcommand{\CoordWt}{\Phi_{0}}
\newcommand{\CoordWtm}{\Phi_{-1}}

\newcommand{\AngAlgGen}{\mathbb{O}}
\newcommand{\AngAlgNormalisedGen}{\hat{\mathbb{O}}}
\newcommand{\RotElementi}{\Theta_i}
\newcommand{\SymGen}{\mathbb{T}}
\newcommand{\NiceBasisSet}{\mathbb{X}}
\newcommand{\NiceNormalisedBasisSet}{{\hat{\mathbb{X}}}}
\newcommand{\VectorSet}{\mathbb{A}}
\newcommand{\TensorNorm}[2]{|#1|_{#2}}
\newcommand{\TensorDNorm}[4]{|#1|_{#2,#3,#4}}
\newcommand{\TensorDSymNorm}[2]{|#1|_{\NiceBasisSet,#2,\SymGen}}
\newcommand{\BigO}[1]{O(#1)}
\newcommand{\ErrorTermsOf}[1]{\BigO{\TensorNorm{#1}{\NiceBasisSet}}}
\newcommand{\ErrorTermskOf}[2]{\BigO{\TensorDNorm{#1}{\NiceBasisSet}{#2}{\SymGen}}}
\newcommand{\ErrorTerms}{\ErrorTermsOf{\MaxF}}
\newcommand{\ErrorTermsk}[1]{\ErrorTermskOf{\MaxF}{#1}}

\newcommand{\GenericNorm}{H[\MaxF](0)}
\newcommand{\normTangN}[1]{\sum_{k=0}^{#1}\GenEnergy{\vecT}[\Lie_\SymGen^k \MaxF](0)}
\newcommand{\normKangN}[1]{\sum_{k=0}^{#1}\GenEnergy{\vecConf}[\Lie_\SymGen^k \MaxF](0)}
\newcommand{\normTunitBifSphN}[1]{\sum_{k=0}^{#1} \GenEnergy{\vecTunit}[\Lie_{\BifurcationCoordinateBasis}^k \MaxF](0)}
\newcommand{\SupInitData}{\sup_{\{0\}\times\Reals^+\times S^2}\sum_i(r^{5/2}\MyFi)^2}

\newcommand{\DivThmRegion}{{\Omega_{(t,\rs)}}}
\newcommand{\DivThmRegionA}{{{\DivThmRegion}_A}}
\newcommand{\DivThmRegionB}{{{\DivThmRegion}_B}}

\newcommand{\SHD}{(*_{S^2})}
\newcommand{\ExtdS}{\Extd_{S^2}}
\newcommand{\LeviCivitaTwoD}{\epsilon}

\newcommand{\dTwo}{d\omega}

\newcommand{\uout}{{u_+}}
\newcommand{\uin}{{u_-}}
\newcommand{\NullInSfc}[1]{{\Sigma^-_{#1}}}
\newcommand{\NullOutSfc}[1]{{\Sigma^+_{#1}}}
\newcommand{\intInSfc}{\int_\NullInSfc{\uout}}
\newcommand{\intOutSfc}{\int_\NullOutSfc{\uin}}
\newcommand{\dThreeUInPOD}{d\uin d^2\omega}
\newcommand{\dThreeUOutPOD}{d\uout d^2\omega}

\newcommand{\SpinTwo}{W}

\newcommand{\PseudoMaxF}{\tilde{F}}
\newcommand{\RiemComp}{R}
\newcommand{\Riem}{\text{Riem}}

\newcommand{\SFn}{u}
\newcommand{\PriceEnergyDensity}{e}
\newcommand{\PriceConfDensity}{e_\mathcal{C}}
\newcommand{\PriceEnergy}{E}
\newcommand{\PriceConf}{E_\mathcal{C}}
\newcommand{\PriceMinEnergy}{E_{\min}}
\newcommand{\SHarmPara}{l}

\newcommand{\LDmult}{\gamma}
\newcommand{\LDwt}{g}
\newcommand{\LDsubwt}{\tilde{g}}
\newcommand{\LDpara}{b}
\newcommand{\LDexp}{\sigma}
\newcommand{\chiLC}{\chi_{\text{LC}}}
\newcommand{\chiLCproto}{\chi_{[-3/4,3/4]}}
\newcommand{\HardyFn}{f}
\newcommand{\sgn}{\text{sgn}}
\newcommand{\chiHardy}{\chi_H}
\newcommand{\PriceLDEnergy}{E_{\LDmult}}

\begin{document}
\maketitle

\begin{abstract}
We study solutions of the decoupled Maxwell equations in the exterior region of a Schwarzschild black hole. In stationary regions, where the Schwarzschild coordinate $r$ ranges over $2M < r_1 < r < r_2$, we obtain a decay rate of $t^{-1}$ for all components of the Maxwell field. We use vector field methods and do not require a spherical harmonic decomposition. 

In outgoing regions, where the Regge-Wheeler tortoise coordinate is large, $\rs>\epsilon t$, we obtain decay for the null components with rates of $|\phi_+| \sim |\alpha| < C r^{-5/2}$, $|\phi_0| \sim |\rho| + |\sigma| < C r^{-2} |t-\rs|^{-1/2}$, and $|\phi_{-1}| \sim |\underline{\alpha}| < C r^{-1} |t-\rs|^{-1}$. Along the event horizon and in ingoing regions, where $\rs<0$, and when $t+\rs>1$, all components (normalized with respect to an ingoing null basis) decay at a rate of $C \uout^{-1}$ with $\uout=t+\rs$ in the exterior region. 
\end{abstract}

\section{Introduction}
\label{sIntro}

The subject of this paper is the study of decay of solutions to the decoupled Maxwell equations in the exterior of a Schwarzschild black hole. The Maxwell field is a $2$-form which we may write in abstract index notation as an antisymmetric $(0,2)$-tensor field on a manifold $\GenericM$, 
\begin{align*}
\MaxF\in&\forms^2(\GenericM) &\text{or}&& \MaxF_{\ia\ib}=&-\MaxF_{\ib\ia} .
\end{align*}
It satisfies the Maxwell equations: 
\begin{align}
*\Extd*\MaxF=&0 &\text{or}&&\CoDeriv^\ia \MaxF_{\ia\ib}=&0 \label{eMaxwellEquationDiv}\\
\Extd\MaxF=&0 && &\CoDeriv_{[\ia} \MaxF_{\ib\ic]}=&0 .\label{eMaxwellEquationAlt}
\end{align}
The exterior region of the Schwarzschild solution is a Lorentz manifold on which the metric is given in terms of coordinates $t\in\Reals$, $r>2M$, $(\theta,\phi)\in S^2$ by
\begin{align}
ds^2 =&-\horifac dt^2 +\horifac^{-1} dr^2 +r^2(d\theta^2 +\sin^2\theta d\phi^2) . 
\label{eSchwarzschildMetric}
\end{align}

This problem comes from general relativity. In general relativity, a model of the universe consists of a space-time manifold $\GenericM$, possibly fields describing matter, and a Lorentz (pseudo-) metric $\g$ which satisfies Einstein's equation. Gravity is described by the curvature of $\g$. The simplest and longest-known solution is Minkowski space, $\Reals^{1+3}$ with the flat metric $-dt^2+dx^2+dy^2+dz^2$. After this, the Schwarzschild manifold is the longest-known solution to Einstein's equation. It is the paradigmatic example of the class of black hole solutions, which play an important role in relativity. The Maxwell field describes electromagnetic radiation. In Einstein's equations, the energy-momentum tensor of the matter fields should influence the curvature. By decoupled, we mean that the electromagnetic field does not influence the Schwarzschild solution, which is taken as a fixed background manifold. We call the Schwarzschild solution the Schwarzschild manifold and use the word solution to refer to solutions of the Maxwell equations \eqref{eMaxwellEquationDiv}-\eqref{eMaxwellEquationAlt}. 

Since $\MaxF$ is a tensor, there is no coordinate independent norm with which to measure it (or, at least, not all components of it). To discuss the decay of $\MaxF$, we make a choice of basis and show that the corresponding components decay. A simple choice of basis consists of the coordinate vector fields rescaled so that they have unit length ($|g(\vecX,\vecX)|=1$). The rescaled vectors are
\begin{align*}
\vecTunit=&\horifac^{-1/2}\dt ,& 
\vecRunit=&\horifac^{1/2}\drOrig ,& 
\vecHunit=&r^{-1}\dtheta ,&
\vecPhunit=&r^{-1}\sin(\theta)^{-1}\dphi.
\end{align*}
Given a time-like vector, there is a natural decomposition of the Maxwell field into electric and magnetic components. Since the Schwarzschild manifold has a time-translation symmetry, this provides a natural choice of time-like direction, $\vecTunit$. The corresponding electric and magnetic components are
\begin{align*}
\Elec_\vecX =& \MaxF_{\vecTunit\vecX} &\vecX\in\{\vecRunit, \vecHunit, \vecPhunit\} ,\\
\Mag_\vecX =& \MaxF_{\vecY\vecZ} & \text{$\vecX, \vecY, \vecZ$ a cyclic permutation of $\vecRunit, \vecHunit, \vecPhunit$} , \\
|\Elec|^2=& |\Elec_\vecRunit|^2 + |\Elec_\vecHunit|^2 + |\Elec_\vecPhunit|^2 ,\\
|\Mag|^2=& |\Mag_\vecRunit |^2 + |\Mag_\vecHunit|^2 + |\Mag_\vecPhunit|^2 .
\end{align*}

Now that we have a choice of components for the Maxwell field, it is possible to state the main decay result of this paper. 
\begin{theorem}[Decay in stationary regions]
\label{tDecayInStationaryRegions}
Let $2M<r_1<r_2<\infty$. There is a constant $C$ and a norm\footnote{The norms used are stated explicitly in section \ref{sStationaryDecay}. For this norm to be finite, it is sufficient that the initial data and its first eight derivatives are bounded and decay like $r^{-(5/2+\epsilon)}$ (see remark \ref{SimplifiedInitialData}). The initial data does not need to decay at the bifurcation sphere, $r\rightarrow2M$. We do not use a spherical harmonic decomposition in our analysis; however, from the structure of the Maxwell equations, spherically symmetric solutions can have no time dependence and cannot decay sufficiently rapidly for the norm $\GenericNorm$ to be finite (see appendix \ref{sExclusionOfNonRadiatable}). } $\GenericNorm$ depending only on $\MaxF$ and its derivatives on the hyper-surface $\{0\}\times(2M,\infty)\times S^2$ such that if $\MaxF$ is a solution to the Maxwell equations \eqref{eMaxwellEquationDiv}-\eqref{eMaxwellEquationAlt}, then for all $t\in\Reals$, $r\in[r_1,r_2]$, $(\theta,\phi)\in S^2$, 
\begin{align*}
|\Elec| +|\Mag|
\leq& C (1+|t|)^{-1}\GenericNorm .
\end{align*}
\end{theorem}

The major advance of this work is to find decay rates which govern all components of the (decoupled) Maxwell field explicitly. The rates we obtain for stationary regions with $r\in(r_1,r_2)$ are significantly slower than the rate of $t^{-5/2}$ which can be obtained in Minkowski space using vector field methods \cite{CK} and the rate of $t^{-3}$ which was derived formally for the Schwarzschild manifold \cite{Price,NewPrice}. Outside of outgoing light-cones, ie where $t<\rs =r +\ln((r-2M)/2M) +C$, decay rates at the same rate as in Minkowski space have already been obtained \cite{IngleseNicolo}. In the outgoing region, we obtain similar results, which we explain below. Certain components of the Maxwell tensor satisfy a scalar wave equation. These components are the zero-weight (spinor or null) components. Previous results for wave equations were sufficiently strong to prove decay for the zero-weight component with a rate of $t^{-1}$ in stationary regions and the appropriate decay in outgoing regions \cite{BlueSterbenz}, although this application was not explicitly stated. $L^\infty_{\text{loc}}$ decay without a rate has also been explicitly obtained using very different techniques \cite{FinsterSmollerMaxwellLG}. The existence and asymptotic completeness of wave operators taking data on the initial surface $t=0$ to the surfaces at $r=2M$ and at infinity has also been shown \cite{Bachelot}. 

Our method starts by using the energy-momentum tensor to generate a positive, conserved energy from the time translation symmetry and a stronger ``conformal energy'' from a vector field $\vecConf$. This follows ideas in \cite{CK} and is very closely related to the analysis of the wave equation in \cite{BlueSofferLongPaper,BlueSterbenz,DafermosRodnianski}. Before the wave estimates were known, a similar method was used \cite{IngleseNicolo}. The growth of the conformal energy is bounded by a ``trapping term'' consisting of the $\Elec_\vecRunit$ and $\Mag_\vecRunit$ components localized near the photon sphere, $r=3M$. In the geometric optics limit, electromagnetic radiation follows null geodesics, which can orbit at $r=3M$. Energy can decay arbitrarily slowly from this region, at least for the wave equation \cite{Ralston}. Thus, it should be expected that there is an obstruction to dispersion near this surface. The trapping term can be controlled because the $\Elec_\vecRunit$ and $\Mag_\vecRunit$ components each satisfy a scalar wave equation of a type that's been previously studied \cite{BlueSterbenz}. This wave equation and the terminology ``zero-weight components'' for $\Elec_\vecRunit$ and $\Mag_\vecRunit$ follow from the analysis of the Price equations \eqref{ePricei}-\eqref{ePriceiv} first appearing in \cite{Price}. We refer to this reduction to scalar wave equations as ``spin reduction''. The control on the conformal energy allows us to conclude:

\begin{lemma}
\label{lSfcDecayInIntro}
There is a constant $C$ and a norm $\GenericNorm$ depending only on $\MaxF$ and its derivatives on the hyper-surface $\{0\}\times(2M,\infty)\times S^2$ such that if $\MaxF$ is a solution to the Maxwell equations \eqref{eMaxwellEquationDiv}-\eqref{eMaxwellEquationAlt}, then for any $2M<r_1<r_2<\infty$ and $t$ sufficiently large, on the surface\footnote{A more general surface is permitted in the statement of lemma \ref{lSurfaceEnergyBound}} $\sfc=\{t\}\times[r_1,r_2]\times S^2$, 
\begin{align*}
\int_{\sfc} \left(|\Elec|^2 +|\Mag|^2\right) \dThreeHF \leq& C t^{-2} \GenericNorm^2 .
\end{align*}
\end{lemma}

From Soolev estimates and integrated decay estimates, like lemma \ref{lSfcDecayInIntro}, for the Lie derivative of $\MaxF$, it is possible to prove pointwise decay estimates. In Minkowski space, the four coordinate directions generate symmetries, so that the Lie derivatives of a Maxwell field also satisfies the Maxwell equations. Although we lack a full set of symmetries, we do have 3 from the time-translation and angular-rotation symmetries. To control a fourth direction, we use the Maxwell equations to ``trade'' the derivatives in the directions of the three symmetries for a radial derivative. With Lie derivatives in all directions controlled, we conclude that theorem \ref{tDecayInStationaryRegions} holds.

\begin{figure}
\begin{center}
\input{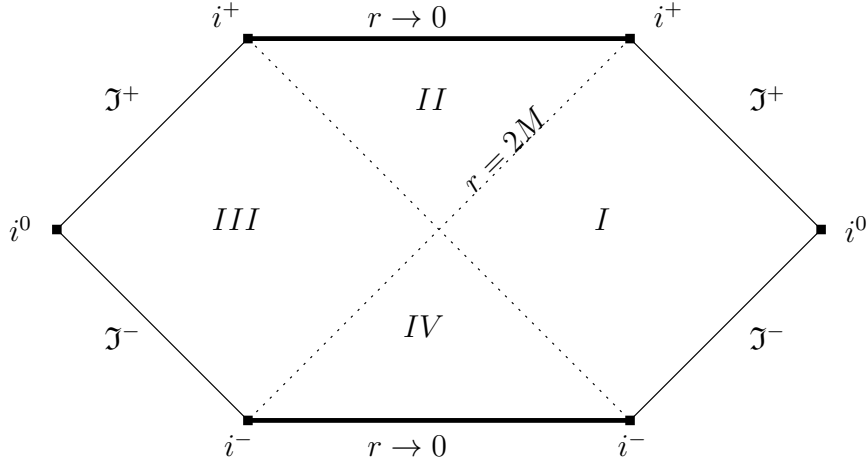}
\caption{A conformal diagram for the maximal extension of the Schwarzschild manifold (suppressing the spherical coordinates). Thin lines represent boundary points at infinity. Thick lines represent the singularity at $r\rightarrow0$. Dotted lines represent the event horizon. Regions $I$ and $III$ are exterior regions, and regions $II$ and $IV$ are interior regions. The surfaces $\scri^\pm$ represent future and past null infinity. The points $i^\pm$ represent future and past null infinity. The points $i^0$ represent spatial infinity. }
\label{figSchwMaxExtension}
\end{center}
\end{figure}

To further explain our results and those of others, we describe the geometry of the Schwarzschild manifold and its importance. This description can be found in most introductory relativity texts (ie \cite{EllisHawking,MTW}). The Lorentz metric is most simply given in terms of coordinates $(t,r,\theta,\phi)$ by \eqref{eSchwarzschildMetric}. As $r\rightarrow\infty$, this metric approaches the flat, Minkowski metric written in spherical coordinates, $ds^2=-dt^2+dr^2+r^2(d\theta^2+sin^2\theta d\phi^2)$. For $r>r_0>2M$, the Schwarzschild solution describes the space-time of a vacuum outside a star of radius $r_0$ and mass $M$. The restriction on $r$ can be relaxed by considering extensions of this manifold. The metric is clearly well-defined in the exterior region $t\in\Reals$, $r\in(2M,\infty)$, $(\theta,\phi)\in S^2$, and in the interior region $t\in\Reals$, $r\in(0,2M)$, $(\theta,\phi)\in S^2$. In the interior region, since $\horifac$ is negative, $r$ is a time-like coordinate, and $t$ is space-like. The maximal analytic extension of any open subset of the Schwarzschild solution is illustrated in the conformal diagram in figure \ref{figSchwMaxExtension}, in which the angular variables are suppressed. There are two exterior regions ($I$ and $III$) and two interior regions ($II$ and $IV$). By an appropriate choice of coordinates, each interior can be smoothly joined to each exterior along a null surface $r=2M$. The manifold is also smooth at the bifurcation sphere where the four regions meet. However, as $r\rightarrow0$, the curvature polynomial $\RiemComp_{\ia\ib\ic\id}\RiemComp^{\ia\ib\ic\id}$ diverges. 

The Schwarzschild manifold is a prototypical solution to Einstein's equation which has inspired many key concepts in general relativity. The asymptotic approach of the metric to the flat, Minkowski metric is known as asymptotic flatness. In the conformal compactification of each exterior region of the Schwarzschild solution, each outgoing geodesic (with $r\rightarrow\infty$ as $t\rightarrow\infty$) ends on future null infinity $\scri^+$, and each ingoing null geodesic (with $r\rightarrow\infty$ as $t\rightarrow-\infty$) starts on past null infinity $\scri^-$. In essence, an asymptotically flat manifold is defined to be one with a future null infinity. A black hole is a region of space-time which cannot be joined by future-directed, null or time-like curves to $\scri^+$, and an event horizon is its boundary. The future, interior region of the Schwarzschild solution is a black hole, and the surfaces where $r=2M$ are the event horizons. In the Schwarzschild manifold, the (future) singularity at $r\rightarrow0$ is separated from $\scri^+$ by the event horizon. The singularity theorems state that under a broad range of conditions, future singularities must form \cite{EllisHawking}. The weak cosmic censorship conjecture asserts that under some genericity condition, which is not yet known, future singularities are always separated from $\scri^+$ by an event horizon. There is a three parameter family of asymptotically flat, known, exact solutions to Einstein's solutions which represent massive, rotating, charged black holes. This is the Kerr-Newman class, and the Schwarzschild solutions are the solutions with positive mass and zero angular momentum and charge. These are stationary, in the sense that they have a time-translation symmetry sufficiently close to null infinity. These solutions also have singularities, but if the angular momentum and charge are beneath a critical threshold, then the singularities are separated from the asymptotically flat regions, in the sense that a future-directed, time-like curve from a point in an exterioir region will either escape to null or space-like infinity or cross the event horizon, but not both. The Kerr-Newman solutions are believed to be the only asymptotically flat, stationary solutions. Physicists believe that all black holes should approach one of the stationary, Kerr-Newman solutions. It is not yet known if a small perturbation of a Cauchy surface for one of the Kerr-Newman solutions will evolve into a solution which remains similar to one of the known solutions. This is the question of black hole stability. 

Stability for Minkowski space was a major and difficult result \cite{CKNL}. Einstein's equations are a complicated system of nonlinear equations in which the geometry is dynamic. The linearization of Einstein's equations about Minkowski space forms a system called the spin 2 field equations. Obtaining decay estimates for the spin $2$ field was one step in this proof \cite{CK}. Decay estimates for the decoupled Maxwell equations were proven at the same time. 

The question of stability of the Schwarzschild solution has also inspired the study of linear fields. In the linearization of Einstein's equation, certain components are determined by the solution to a simple wave equation \cite{ReggeWheeler}, and the remaining components are determined by the solution to a more complicated wave equation \cite{Zerilli}. Using spinors, Price was able to present a more unified presentation for all components of several important, physical systems, according to their spin. Any wave equation is said to have spin $0$. The Dirac system has spin $1/2$. The Maxwell field has spin $1$. For any solution of Einstein's equations in vacuum, the non-vanishing components of the curvature satisfy certain relations from the Bianchi identities and Einstein's equation. In Minkowski space, since the curvature is zero, solutions to the linearization of Einstein's equation satisfy the same relations, which are called the spin $2$ field equations. In the Schwarzschild manifold, since the curvature is non-vanishing, in the linearization of Einstein's equations, there are additional terms arising from the derivative with respect to the perturbed metric of the original Christoffel terms. Thus, Price distinguishes between the spin $2$ field equations and the linearization of Einstein's equation. Formal arguments suggest a rate of $t^{-3}$ for fields of all integer spin and for the linearization of Einstein's equations \cite{Price}. A similar, spinorial presentation of these systems has been made for the Kerr-Newman solutions \cite{Teukolsky}. In each case, certain components were found to satisfy scalar wave equations and then acted as forcing terms in the equations governing the remaining components. 

Most of the subsequent analysis has been focused on the decoupled wave equation. The literature is vast, and we list only some of the results. Solutions are known to remain uniformly bounded in time \cite{KayWald}. The scattering theory, concerning the map from the initial data to the the limit on $\scri^+$ and the event horizon, has also been studied on the Schwarzschild manifold \cite{DimockKayScattering} and on the more general Kerr-Newman solutions \cite{Hafner}. On the Schwarzschild manifold, vector-field techniques have been used to obtain decay results in three main steps \cite{BlueSterbenz,BlueSofferLongPaper,DafermosRodnianski}. First, the vector field $\vecConf$ is used to introduce a conformal energy, which is not conserved because of trapping. Second, a radial vector field is used to prove a local decay estimate to control the trapping term. In $\Reals^{1+3}$, the radial derivative can be used to make a somewhat similar estimate \cite{Morawetz}. In $\Reals^{1+3}$, estimates involving $\vecConf$ and the radial vector field are both referred to as Morawetz estimates. In this step, a spherical harmonic decomposition was used in the proofs, but this is no longer necessary \cite{DafermosRodnianskiNoSHDecomp}. Because the scalar equation governing the zero-weight components of the Maxwell field has a simple structure, in appendix \ref{sOneDWaveAnalysis}, we are able to modify the earlier method to obtain decay without using a spherical harmonic. Third, the conformal energy is used to control norms. In \cite{BlueSterbenz,DafermosRodnianski}, a strong local decay estimate is proven and additional angular derivatives are used to obtain a $L^\infty_{\text{loc}}$ decay rate of $t^{-1}$ and a similar decay rate in outgoing regions. In \cite{DafermosRodnianski}, an additional vector field, $Y$, is used to also prove decay estimates along the event horizon. These require weighted $H^4$ or $H^5$ norms of the initial data to be bounded. By an $H^k$ norm, we mean, roughly speaking, that the $k$th derivative of a solution $u$ is square integrable. In \cite{BlueSofferLongPaper}, only weighted $H^{1+\epsilon}$ norms of the initial data are needed, but a weaker local decay estimate was obtained, which led to less control on the conformal energy and a decay rate of $t^{-1/3}$ for a weighted, spatial $L^6$ norm. Using an entirely different technique, based on a representation of the propagator, $L^\infty_\text{loc}$ decay has been proven for the wave equation on subcritical Kerr-Newman solutions \cite{FinsterKamranSmollerYauWave}. 

The spin $1/2$ system is the Dirac model for the electron. On the sub-critical Kerr-Newman solutions, scattering results \cite{HafnerNicolasDirac} and $L^\infty_{\text{loc}}$ decay \cite{FinsterKamranSmollerYauDirac} have also been proven. 

For the linearization of Einstein's equations about the Kerr-Newman solutions, the equations found in \cite{Teukolsky} were found to have no unstable modes \cite{Whiting}. For the linearization about the Schwarzschild solution, the simpler equations in \cite{ReggeWheeler} satisfy an integrated decay estimate \cite{BlueSofferReggeWheeler} and $L^\infty_{\text{loc}}$ decay \cite{FinsterSmollerMaxwellLG}. Although the application is not explicitly stated, decay at a rate of $t^{-1}$ follows from \cite{BlueSterbenz}. 

For the full Einstein equations on a black hole background, decay results are known in the spherically symmetric case, when the Einstein equations are coupled to a scalar wave equation and the Maxwell field \cite{DafermosRodnianskiSpherical}. By Birkhoff's theorem (see \cite{EllisHawking}), the Schwarzschild manifold is the only spherically symmetric solution to Einstein's equation (treating Minkowski space as the special sub-case with $M=0$). The decay rate obtained for the scalar field is $\uout^{-3+\epsilon}$ along the event horizon. The $\uout^{-3+\epsilon}$ decay rate is known as Price's law and has important implications for the strong cosmic censorship conjecture \cite{DafermosCensorship}. 

To discuss decay outside regions of fixed $r$, it is necessary to introduce components with respect to a null tetrad, a basis built from null vectors with certain properties. Physicists may know these as spinor components \cite{Price,Stewart}, and mathematicians, as the null decomposition \cite{CK}. We present one tetrad here and discuss exactly what we mean by a null tetrad in section \ref{sNotation}. 

We start by introducing the Regge-Wheeler radial coordinate, $\rs$, defined by
\begin{align*}
\frac{d r}{d\rs}=& \horifac ,
& r(0)=& 3M .
\end{align*}
The exterior region of the Schwarzschild solution is given by $(t,\rs,\theta,\phi)$ ranging over $\Reals\times\Reals\times S^2$. In these coordinates, the Lorentz metric becomes
\begin{align*}
\g=&-\horifac dt^2 +\horifac d\rs^2 +r^2(d\theta^2 +\sin^2\theta d\phi^2) . 
\end{align*}
From this form of the metric, it is clear that any multiple of the vectors $\dt\pm\dr$ are null. To define our null tetrad, we use $\vecea$ and $\veceb$ to denote an orthonormal basis of tangents vectors to $S^2$ and $\LeviCivitaTwoD$ to denote the antisymmetric, Levi-Civita tensor on $S^2$. The null tetrad we will use to state our results is
\begin{align*}
&(\dt+\dr) ,&
&\horifac^{-1}(\dt-\dr) ,&
&r^{-1} \vecea , &
&r^{-1} \veceb .
\end{align*}
The covariant derivative of this null tetrad along ingoing, radial, null geodesics is zero. Thus, having found a natural choice of null tetrad on the initial surface $t=0$, we have extended it to the entire future of the initial surface by parallel transport along null geodesics falling into the black hole. This is useful for considering limits as $r\rightarrow 2M$. Had we extended the basis by parallel transport along outgoing null geodesics, to study the problem as $r\rightarrow\infty$, the factor of $\horifac^{-1}$ would have been on $(\dt+\dr)$ instead of $(\dt-\dr)$. However, since $\horifac\rightarrow1$ as $r\rightarrow\infty$, the difference between our choice of null tetrad and the natural choice is vanishingly small. Therefore, we use our choice of null tetrad throughout the future, $t\geq0$. 

The decay of these spinor or null components is not simply a decay in time. This is known from the behavior in $\Reals^{1+3}$. In that case, the heuristic is that the bulk of solutions to the Maxwell equations travel out along the light-cone $t\sim|\vec{x}|$. In any fixed region, there is decay because the wave leaves the region. As the light-cone expands, the average value of the wave intensity drops. Moving with the wave, the intensity decays as it is spread over the increasing area of the light-cone. Thus, some of the decay occurs as a result of the wave being far from the light cone, and some occurs as a result of the light-cone being very large. Similar behavior occurs on the Schwarzschild manifold. The null coordinates
\begin{align*}
\uout=& t+\rs &\text{and}&& 
\uin=& t-\rs 
\end{align*}
are used to measure the distance from the light-cone, and $r$ is used to measure the size of the light-cone in the outgoing direction. In the ingoing direction, the radius of the surface of the light-cone also goes like $r$, but since this approaches $2M$, the decay occurs only in the null coordinates. 

\begin{theorem}[Decay outside stationary regions]
\label{tNearAndFarDecay}
There is a constant $C$ and a norm $\GenericNorm$ depending only on $\MaxF$ and its derivatives on the hyper-surface $\{0\}\times(2M,\infty)\times S^2$ such that if $\MaxF$ is a solution to the Maxwell equations \eqref{eMaxwellEquationDiv}-\eqref{eMaxwellEquationAlt}, then for all $t\geq0$, $\rs>1$, $(\theta,\phi)\in S^2$, 
\begin{align*}
|\MaxF(\dt+\dr,r^{-1}\vecea)|\leq& C r^{-3/2}\uout^{-1} \GenericNorm ,\\
|\MaxF(\dt+\dr,\horifac^{-1}(\dt-\dr))|+|\MaxF(r^{-1}\vecea,r^{-1}\veceb)\LeviCivitaTwoD^{\isa\isb}|  \leq& C r^{-2}\left(\frac{\uout-|\uin|}{\uout(1+|\uin|)}\right)^{1/2} \GenericNorm ,\\
|\MaxF(\horifac^{-1}(\dt-\dr),r^{-1}\vecea)| \leq& Cr^{-1}(1+|\uin|)^{-1} \GenericNorm .
\end{align*}
If we restrict to $\uin\leq 0$, then 
\begin{align*}
|\MaxF(\horifac^{-1}(\dt-\dr),r^{-1}\vecea)| \leq& Cr^{-1}(1+|\uin|)^{-3/2} \GenericNorm .
\end{align*}
Under the same hypotheses, then for all $t\geq0$, $\rs<-1$, $(\theta,\phi)\in S^2$ such that $\uout>1$, 
\begin{align*}
|\MaxF(\dt+\dr,r^{-1}\vecea)| \leq& C \uout^{-1}\GenericNorm ,\\
|\MaxF(\dt+\dr,\horifac^{-1}(\dt-\dr))|+|\MaxF(r^{-1}\vecea,r^{-1}\veceb)\LeviCivitaTwoD^{\isa\isb}| \leq& C \uout^{-1}\GenericNorm ,\\
|\MaxF(\horifac^{-1}(\dt-\dr),r^{-1}\vecea)| \leq& C \uout^{-1}\GenericNorm .
\end{align*}
\end{theorem}

This gives a $\uout^{-1}$ decay rate for all components (since, either $r>C\uout$ or $\uin>C\uout$ in the far region $\rs>1$). Outside the outgoing light-cone, where $0<t<\rs-1$, the decay rates are $r^{-5/2}$, $r^{-2}\uin^{-1/2}$, and $r^{-1}\uin^{-3/2}$. These are the same rates as can be obtained in Minkowski space using vector field methods, with $\uin=t-|\vec{x}|$ in Minkowski space. (Faster decay rates can be obtained using conformal compactification and other methods.) Approaching null infinity inside the light-cone, with $(1+\epsilon)\rs>t>\rs>0$, the decay rates are $r^{-3/2}\uout^{-1}$, $r^{-2}\uin^{-1/2}$, and $r^{-1}\uin^{-1}$. Thus, the decay rates for the first two components are the same as in Minkowski space, but the last component decays more slowly than in Minkowski space. The slow decay for this component comes from the slow decay rate of $t^{-1}$ in stationary regions.  

In addition to the Maxwell equations, one can imagine studying the spin $2$ field equations. A spin $2$ field is a $(0,4)$ tensor with the following symmetries 
\begin{align}
\SpinTwo_{\ib\ia\ic\id}=&-\SpinTwo_{\ia\ib\ic\id} \label{eSpinTwoDefnFirst}\\
\SpinTwo_{\ia\ib\id\ic}=&-\SpinTwo_{\ia\ib\ic\id} \label{eSpinTwoASSecond}\\
\SpinTwo_{[\ia\ib\ic]\id}=&0 \label{eFirstBianchiSpinTwo}\\
\SpinTwo_{\ia\ib\ic}{}^\ia=& 0 \label{eSpinTwoTraceFree} ,
\end{align}
and which satisfies the spin $2$ field equations
\begin{align}
\CoDeriv^{\ic}\SpinTwo_{\ic\id\ia\ib}=&0  \label{eSpinTwoDivFree} \\
\CoDeriv_{[\ie}\SpinTwo_{\ic\id]\ia\ib}=& 0 .\label{eSpinTwoDefnLast} 
\end{align}
The symmetries of a spin $2$ field are similar to the antisymmetry of a Maxwell field, and the spin $2$ field equations are similar to the Maxwell equations. If the vacuum Einstein equations are satisfied, then the Ricci curvature vanishes, and the Weyl curvature satisfies the spin $2$ field equations. In $\Reals^{1+3}$, the spin $2$ field equations are a good model for the linearization of Einstein's equation about the Minkowski solution, but this is not true for the linearization about other solutions. In Cartesian coordinates on $\Reals^{1+3}$, the Christoffel symbols and the curvature are zero. If one introduces a perturbed metric on Minkowski space and treats the Weyl tensor as a tensor field on the original space-time, then the difference between the covariant derivative of the Weyl tensor with respect to the perturbed metric and the original metric will be second order in the perturbation. Thus, ignoring second order terms, the perturbed Weyl tensor satisfies the spin $2$ field equations on the original metric. In this sense, the spin $2$ field equations are the linearization of the vacuum-Einstein equation about Minkowski space. This is the motivation for studying the spin $2$ field in \cite{CK}. When linearizing around a curved space-time, the Christoffel symbols do not vanish, and the linearized Einstein equations do not reduce to the spin $2$ field equations. More drastically, there is a Buchdahl constraint \cite{Stewart} from applying two covariant derivatives, two contractions, and the spin $2$ field equations, 
\begin{align*}
\Riem^{\ic\id\ie}{}_{(\ia} \SpinTwo_{\ib)\ie\ic\id} = 0 ,
\end{align*}
where $\Riem$ is the Riemann curvature of the background. On the Schwarzschild manifold, this forces $\SpinTwo(\vecTunit+\vecRunit,\vecTunit-\vecRunit,\vecTunit\pm\vecRunit,r^{-1}\vecea)$ to vanish everywhere. This forces the derivative of other components to vanish, so that there is only a finite dimensional family of spherically symmetric solution. These are similar to the spherically, and non-decaying solutions, which are discussed in appendix \ref{sExclusionOfNonRadiatable}. 

Nonetheless, one can ignore the Buchdahl constraint and study the spin $2$ field equations. The system has been studied formally as a system of transport equations \cite{Price}. Here we report that it is possible to use a method similar to that we used for the Maxwell equations. It is well-known in the literature that the Bel-Robinson can be used to define a conserved, positive-definite quantity from a time-like Killing vector in the same way that the energy-momentum tensor can for the Maxwell field. In addition to the conserved energy generated this way, one can use a quantity defined in terms of the time-translation symmetry and the vector field $\vecConf$. One can again use a method of ``spin-reduction'' to introduce a ``pseudo-Maxwell tensor'', $\PseudoMaxF_{\ia\ib}=\SpinTwo_{\ia\ib\ic\id}\vecTunit^\ic\vecRunit^\id$, which satisfies the Maxwell equations. One can use control of the pseudo-Maxwell tensor to control the trapping terms for the spin $2$ field to control integrated norms of the spin $2$ field. One can then use the symmetries of the Schwarzschild space-time, the field equations, and Sobolev estimates to prove $t^{-1}$ pointwise-in-time decay for the non-spherically symmetric components of the spin $2$ field. Clearly this is pointless, since the dynamics of the spin $2$ field are trivial. However, we expect that a similar analysis will apply to the genuine, linearized gravity system. The linearized gravity equations are more complicated than the spin $2$ field equations because there are terms involving the perturbed Christoffel symbols contracted against the unperturbed and nonvanishing Weyl tensor.

In section \ref{sNotation}, we introduce several sets of vector fields to provide a simpler notation for discussing the null decomposition of the Maxwell field and symmetries. The null decomposition and spinor decomposition are essentially equivalent. We estimate energies in section \ref{sEnergies}. We first review the use of the energy-momentum tensor, and then use it to define an energy and a weighted, conformal energy. The growth of the conformal energy is controlled by a trapping term which depends only on the zero-weight component. These satisfy a simple wave equation, which allows us to bound the conformal energy. In section \ref{sStationaryDecay}, we use this bound and trade Lie derivatives to prove theorem \ref{tDecayInStationaryRegions}. In section \ref{sMovingDecay}, we use the bounds and integration along null geodesics to prove theorem \ref{tNearAndFarDecay}. In appendix \ref{sExclusionOfNonRadiatable}, we show that there are no spherically symmetric components of the Maxwell field which decay sufficiently rapidly at infinity. In appendix \ref{sOneDWaveAnalysis}, we analyse the wave equation governing the zero-weight component using a simplified version of the arguments in \cite{BlueSterbenz}. This simplified version does not require a spherical harmonic decomposition.

%
%

\section{Notation}
\label{sNotation}
The main purpose of this section is to collect various vector fields and components of the Maxwell field, so that the reader can compare the notation used in different places in this paper and elsewhere. 

We begin with some simple notation. We sometimes use vectors as indices on tensors to denote the corresponding component. We use the notation $\g$ to denote the metric, $\omega=(\theta,\phi)\in S^2$ for the angular coordinate, $\dAng$ for the angular derivative, and $\Omega$ to denote the metric on $S^2$. Except in the statement of results, we use $(t,\rs,\theta,\phi)$ coordinates, unless otherwise specified. We use $\formXg$ to denote the one form generated by lowering a vector field $\vecX$ with the metric. We use the insertion operation $\insertion{\vecX}$ which takes the $(0,m)$ tensor $\tensorA$ to the $(0,m-1)$ tensor $\insertion{\vecX}\tensorA$ defined by $\insertion{\vecX}\tensorA(\vecY_1,\ldots,\vecY_{m-1})=\tensorA(\vecX,\vecY_1,\ldots,\vecY_{m-1})$. We use $C$ to denote an arbitrary constant which may change from line to line in a calculation. 

We will prove estimates for $t\geq0$. Because the Schwarzschild manifold is time symmetric, similar estimates hold for $t\in\Reals$. In particular, estimates in stationary regions will remain the same, and in ingoing and outgoing directions, $t$ and $\dt$ must be replaced by $-t$ and $-\dt$.

\subsection{Coordinates, Bases, and Field Components}
\label{ssComponents}
Recall the definition of the coordinates $t$, $\rs$, $\theta$, $\phi$, $\uin$, and $\uout$ from the introduction. 

We frequently use the coordinate vector fields
\begin{align*}
\vecT=&\dt ,& 
\vecR=&\dr ,& 
\vecH=&\dtheta ,&
\vecPh=&\dphi, 
\end{align*}
and the corresponding normalized vector fields
\begin{align*}
\vecTunit=&\horifac^{-1/2}\dt ,& 
\vecRunit=&\horifac^{-1/2}\dr ,& 
\vecHunit=&r^{-1}\dtheta ,&
\vecPhunit=&r^{-1}\sin(\theta)^{-1}\dphi .
\end{align*}
From the definition of $\rs$, the definition of $\vecRunit$ given here is the same as the one given in the introduction. 

We also use null tetrads. In the standard presentation of a null tetrad, the tangent space is complexified. A basis $\{\vecl, \vecn, \vecm, \vecmb\}$ in which $\vecl$ and $\vecn$ are (real) null vectors, $\vecmb$ is the complex conjugate of $\vecm$, $\g(\vecl,\vecn)=-2$, $\g(\vecm,\vecmb)=2$, and all other inner products between the basis vectors are zero. If $\vecX$ and $\vecY$ are unit vectors orthogonal to $\vecl$ and $\vecn$, a null tetrad can be defined by taking $\vecm = \vecX+ i\vecY$. Because of this, we will also call a basis $\{\vecl, \vecn, \vecX, \vecY\}$ a null tetrad if $\g(\vecl,\vecn)=-2$, $\g(\vecX,\vecX)=\g(\vecY,\vecY)=1$, and all other inner products are zero. We will generally ignore the distinction between the two definitions of a null tetrad. 

One advantage of null tetrads is that they assign weights to certain quantities. Rescaling $\vecl$ and $\vecn$ by $\lambda$ and $\lambda^{-1}$ respectively preserves the null tetrad structure, as does rescaling $\vecm$ by $e^{is}$ (and $\vecmb$ by the conjugate, $e^{-is}$). If, under such a change of basis, a quantity transforms as a power of $\lambda$ or of $e^{-is}$ then the corresponding powers are the conformal and spin weights of the quantity. 

We use several null tetrads. The first is the ``stationary'' tetrad: 
\begin{align*}
\veclunit=&\vecTunit+\vecRunit = \horifac^{-1/2}(\dt +\dr) \\
\vecnunit=&\vecTunit-\vecRunit = \horifac^{-1/2}(\dt -\dr) \\
\vecm=& \vecHunit +i\vecPhunit = \frac1r\dtheta +\frac{i}{r\sin\theta}\dphi  \\
\vecmb=& \vecHunit -i\vecPhunit =\frac1r\dtheta -\frac{i}{r\sin\theta}\dphi  .
\end{align*}
Price \cite{Price} uses a basis which is parallelly transported along outgoing, radial, null geodesics, $\gamma(s)=(s,s+C,\theta_0,\phi_0)$ (in the $(t,\rs,\theta,\phi)$ coordinates). The null vectors are  
\begin{align*}
\veclP=& \horifac^{-1}(\vecT +\vecR) , \\
\vecnP=& \vecT -\vecR , 
\end{align*}
and the angular basis vectors remain the same. To prove theorem \ref{tNearAndFarDecay}, in section \ref{sMovingDecay}, we use a basis adapted to ingoing, radial, null geodesics. Certain expressions are simplified by using the following coordinate-like vector fields. 
\begin{align*}
\veclcoord=& \dt+\dr =\vecT +\vecR,\\
\vecncoord=& \dt-\dr=\vecT -\vecR ,\\
\vecmcoord=& \vecH +\frac{i}{\sin\theta} \vecPh  =\dtheta +\frac{i}{\sin\theta}\dphi .
\end{align*}
Christodoulou and Klainerman \cite{CK} avoid complexifying the tangent space by using an orthonormal basis tangent to the sphere at each point. We use $\vecea$ and $\veceb$ to denote an orthonormal basis on $S^2$. Thus, $r^{-1}\vecea$ and $r^{-1}\veceb$ are unit vectors in the Schwarzschild manifold. The indices $\isa,\isb,\ldots$ are used for directions tangent to the sphere. In summary, we have three null tetrads and a coordinate null basis, 
\begin{align}
&\{\veclunit,\vecnunit,\vecm,\vecmb\} ,\label{eStationaryTetrad}\\
&\{\veclP,\vecnP,\vecm,\vecmb\} ,\label{ePriceTetrad} \\
&\{\veclunit,\vecnunit,r^{-1}\vecea,r^{-1}\veceb\} ,\label{eCKTetrad} \\
&\{\veclcoord,\vecncoord,\vecmcoord,\vecmbcoord\} .\label{eCoordTetrad}
\end{align}

The bases can be used to define the corresponding components of the Maxwell field. The electric and magnetic decomposition was already explained in the introduction. We now introduce a null decomposition and spinor components. These are very closely related but differ in the notation and slightly in the definition. The null decomposition consists of two scalars, $\rho$ and $\sigma$, and two $1$-forms tangent to spheres, $\CKalpha$ and $\CKalphab$. The spinor components are three complex-valued functions. These are defined in terms of the tetrad in \eqref{eCKTetrad} and in \eqref{eStationaryTetrad} by
\begin{align*}
\CKalpha(\vecea) =&\MaxF(\veclunit,\vecea) &
\MyFp =& \MaxF(\veclunit,\vecm) \\
\CKrho   =& \frac12 \MaxF(\veclunit,\vecnunit) &
\MyF =&\frac12 (\MaxF(\veclunit,\vecnunit) +i\MaxF(\vecmb,\vecm))\\
\CKsigma =& \frac12 \MaxF(\vecea,\veceb)\LeviCivitaTwoD^{\isa\isb} &
\\
\CKalphab(\vecea)=&\MaxF(\vecnunit,\vecea) ,&
\MyFm=&\MaxF(\vecnunit,\vecmb) .
\end{align*}
The spin component index in $\MyFi$ refers to both the conformal and spin weight. These components are related by
\begin{align*}
\MyFp=& \CKalpha(\vecm)     & |\MyFp|^2=&|\CKalpha|^2 \\
\MyF =& \CKrho +i\CKsigma   & |\MyF|^2 =&|\CKrho|^2+|\CKsigma|^2 \\
\MyFm=& \CKalphab(\vecmb)   & |\MyFm|^2 =&|\CKalphab|^2 .
\end{align*}
The null decomposition, $\CKalpha$, $\CKrho$, $\CKsigma$, and $\CKalphab$, more accurately represents the geometric behavior of the components. The spinor notation reveals the spin and conformal weight more easily, simplifies several expressions, and suggests connections between the spin $0$ wave equation, the spin $1$ Maxwell equation, and the spin $2$ equations. We typically write expressions in terms of the spinor components but think in terms of the null decomposition. 

The spinor components in \cite{Price} are slightly different from the ones we use. Since $\MyFi$ has conformal weight $i$, replacing the null tetrad \eqref{eStationaryTetrad} by \eqref{ePriceTetrad} will take $\MyFi$ to $\horifac^{i/2}\MyFi$. These are the components initially used in \cite{Price}. 

The spinor components are related to the electric and magnetic components by 
\begin{align*}
\MyFp
=&(\Elec_\vecHunit+\Mag_\vecPhunit) + i(\Elec_\vecPhunit-\Mag_\vecHunit)\\
\MyF 
=&\Elec_\vecRunit +i\Mag_\vecRunit \\
\MyFm
=&(\Elec_\vecHunit-\Mag_\vecPhunit) - i(\Elec_\vecPhunit+\Mag_\vecHunit) .\\
|\MyFp|^2 +2|\MyF|^2 + |\MyFm|^2
=&|\CKalpha|^2+2(|\CKrho^2|+|\CKsigma|^2)+|\CKalphab|^2  
= 2(|\Elec|^2+|\Mag|^2) .
\end{align*} 

Certain calculations are simplified by using the null basis \eqref{eCoordTetrad}. We define Maxwell field components associated to this null basis by
\begin{align*}
\CoordWtp =& \MaxF(\veclcoord,\vecmcoord) 
&=&r\horifac^{1/2}\MyFp \\
\CoordWt  =& \frac12(\MaxF(\veclcoord,\vecncoord)\horifac^{-1}r^2 +\MaxF(\vecmbcoord,\vecmcoord) )
&=&r^2\MyF\\
\CoordWtm =& \MaxF(\vecncoord,\vecmbcoord) 
&=&r\horifac^{1/2}\MyFm .
\end{align*}

\newcommand{\Buout}{U_+}
\newcommand{\Buin}{U_-}
To discuss the maximally extended Schwarzschild solution in a neighborhood of the bifurcation sphere, it is typical to introduce coordinates
\begin{align*}
\Buout=& e^{\uout/4M} ,\\
\Buin=& -e^{-\uin/4M} .
\end{align*}
In the exterior region, these range over $\Buout\in(0,\infty)$ and $\Buin\in(-\infty,0)$. The coordinates $(\Buout,\Buin,\theta,\phi)$ can be used in a neighborhood of the bifurcation sphere, and the bifurcation sphere corresponds to $(\Buout,\Buin)=(0,0)$. The surface $t=0$ corresponds to $\Buout\Buin=-1$. This surface extends through the bifurcation sphere to the surface $t=0$ in the other exterior  region. Since
\begin{align*}
e^\frac{\uout+\uin}{4M}=& C r e^{\frac{r}{2M}}\horifac , 
\end{align*}
on the initial surface $t=0$, the coordinate vector fields are
\begin{align*}
\frac{\partial}{\partial \Buout}=& C r^{1/2} e^\frac{r}{4M} \horifac^{1/2} (\vecT +\vecR ), \\
\frac{\partial}{\partial \Buin}=& C r^{1/2} e^\frac{r}{4M} \horifac^{1/2} (\vecT -\vecR ). 
\end{align*}
Thus, on the initial data surface and near the bifurcation sphere, the coordinate vector fields $\frac{\partial}{\partial \Buout}$ and $\frac{\partial}{\partial \Buin}$ are related to $\veclunit$ and $\vecnunit$ by bounded, nonvanishing functions. If $\{\frac{\partial}{\partial \Buout},\frac{\partial}{\partial \Buin},\vecH,\vecPh\}$ are used to define a tetrad, the corresponding components of the Maxwell field are equivalent to $\MyFi$. Since these vector fields are coordinate vector fields, they commute. To restrict attention to the region near the bifurcation sphere, we will often apply smooth, cut-off functions $\chi_{<0}(\uout)$ and $\chi_{<0}(\uin)$ which are smooth, identically zero for $\uout>1$ and $\uin>1$ respectively, and identically one for $\uout>0$ and $\uin>0$ respectively. The vector fields $\chi_{<0}(\uout)\frac{\partial}{\partial \Buout}, \chi_{<0}(\uin)\frac{\partial}{\partial \Buin},\vecH,\vecPh$ still commute.

\newcommand{\BifurcationCoordinateBasis}{\hat{\tilde{\mathbb{X}}}}

\subsection{Norms and Lie derivatives}
\label{ssLieDerivatives}
With the goal of applying derivatives to the components of the Maxwell tensor, we introduce several collections of vector fields. Since the vector fields $\vecH$ and $\vecPh$ are not smooth, we use the three rotations of $S^2$ about the coordinate axes, $\RotElementi$. We treat these as vector fields on the Schwarzschild manifold. The collections of vector fields we will use are
\begin{align*}
\AngAlgGen              =& \{ \RotElementi \} ,\\
\AngAlgNormalisedGen    =& \{ r^{-1} \RotElementi \} ,\\
\SymGen                 =& \{ \vecT,\RotElementi \} ,\\
\NiceBasisSet           =& \{ \vecR, \vecT,\RotElementi \} ,\\
\NiceNormalisedBasisSet =& \{ \vecRunit, \vecTunit, r^{-1} \RotElementi \} , \\
\BifurcationCoordinateBasis =& \{ \chi_{<0}(\uout)\frac{\partial}{\partial \Buout}, \chi_{<0}(\uin)\frac{\partial}{\partial \Buin},\RotElementi \} .
\end{align*}
Since the Schwarzschild manifold is static and spherically symmetric, $\SymGen$ generates symmetries of the space-time. The normalized vectors in $\NiceNormalisedBasisSet$ are used to define the norms of the electric and magnetic components of the Maxwell tensor. (The three $\RotElementi$'s can be used to define corresponding components of the electric and magnetic field. Taken together these give the angular components.) On the initial data hypersurface $t=0$ and near the bifurcation sphere (ie, where $\rs<0$), the coordinate vectors in $\BifurcationCoordinateBasis$ can also be used to define the Maxwell field components. 

We now recall some convenient notation for discussing collections of vectors and scalar functions from \cite{CK}. For two sets of vector fields, $\VectorSet_i$, the covariant and Lie derivatives are 
\begin{align*}
\Lie_{\VectorSet_1} \VectorSet_2 =& \{ \Lie_{\vecX_1} \vecX_2 | \vecX_i\in\VectorSet_i \} ,\\
\CoDeriv_{\VectorSet_1} \VectorSet_2 =& \{ \CoDeriv_{\vecX_1} \vecX_2 | \vecX_i\in\VectorSet_i \} .
\end{align*}
For two such sets and a $(0,2)$ tensor $\tensorA$, the components of $\tensorA$ with respect to the vector fields are the collection of scalar functions
\begin{align*}
\tensorA(\VectorSet_1,\VectorSet_2)=& \{ \tensorA(\vecX_1,\vecX_2) | \vecX_i\in\VectorSet_i \} .
\end{align*}
Similarly, for a set of vectors $\VectorSet$ and a collection of scalar functions $\{f_i \}$, the derivatives are defined as 
\begin{align*}
\Lie_\VectorSet\{f_i\}
=\CoDeriv_\VectorSet\{f_i\}
= \VectorSet \{f_i\}
=& \{ \vecX f | \vecX\in\VectorSet, f\in\{f_i\} \}. 
\end{align*}
This definition holds since the Lie, covariant, and directional derivatives are the same operation when applied to scalar functions. For tensor fields, a similar notation can be used to generate collections of tensor fields and to consider their components. For example, 
\begin{align*}
(\Lie_{\VectorSet_1} \tensorA)(\VectorSet_2,\VectorSet_3)
=& \{ (\Lie_{\vecX_1} \tensorA)(\vecX_2,\vecX_3) |\vecX_i\in\VectorSet_i \} .
\end{align*}
The same can be defined for iterated Lie or covariant derivatives. 

The norm of a $1$-form or a $(0,2)$ tensor with respect to a set of vector fields is
\begin{align*}
\TensorNorm{\formZg}{\VectorSet}
=& \sum_{\vecX\in\VectorSet} |\formZg(\vecX)| , \\
\TensorNorm{\tensorA}{\VectorSet}
=& \sum_{\vecX,\vecY\in\VectorSet} |\tensorA(\vecX,\vecY)| .
\end{align*}
The $n$-derivative norms of a $(0,m)$ tensor with respect to components in $\VectorSet_1$ and derivatives in the $\VectorSet_2$ directions are defined to be
\begin{align*}
\TensorDNorm{\tensorA}{\VectorSet_1}{n}{\VectorSet_2}^2
=& \sum_{k=0}^{n} \TensorNorm{ \Lie_{\VectorSet_2}^k \tensorA} {\VectorSet_1}^2 \\
=& \sum_{k=0}^{n} \sum_{\vecX_1,\ldots,\vecX_k\in\VectorSet_2, \vecY_1,\ldots,\vecY_m\in\VectorSet_1} | (\Lie_{\vecX_k}\ldots\Lie_{\vecX_1}\tensorA)(\vecY_1,\ldots,\vecY_m) |^2 .
\end{align*}
We note that this notation can be applied equally well with $S^2$ tangent $1$-forms, such as $\CKalpha$ and $\CKalphab$ as any other forms.

%
%
\section{Control of Energies}
\label{sEnergies}
\subsection{Summary of the Lagrangian method}
We briefly outline the Lagrangian theory for a general field $\GeneralField$ on a manifold, $\GenManifold$. It is assumed that there is a scalar Lagrangian $\Lagrangian[x,\phi,\CoDeriv\phi]$ which is used to define the action, 
\begin{align*}
\Action=& \int_{\GenManifold} \Lagrangian[x,\phi,\nabla\phi] \dFourGen.
\end{align*}
If $\GeneralField$ is a minimizer (or, more generally, a critical point) of the action, then $\GeneralField$ will satisfy the Euler-Lagrange equation
\begin{align*}
\frac{\VarDeriv\Lagrangian}{\VarDeriv\GeneralField}-\CoDeriv^\ia\frac{\VarDeriv\Lagrangian}{\VarDeriv \CoDeriv^\ia\GeneralField} =&0 .
\end{align*}
One can then define the energy-momentum tensor from this
\begin{align*}
\StressEnergy_{\ia\ib}
=& \frac12 \left(\CoDeriv_{\ia}\GeneralField \frac{\VarDeriv\Lagrangian}{\VarDeriv \CoDeriv^\ia\GeneralField} -\g_{\ia\ib}\Lagrangian \right) ,
\end{align*}
which, by the Euler-Lagrange equation, satisfies
\begin{align}
\CoDeriv^{\ia}\StressEnergy_{\ia\ib}=&0 .
\label{eEMSDivFree}
\end{align}
For any vector field $\vecX$, the generalized momentum vector $\GenMomentum{\vecX}$ and deformation $2$-tensor $\Deformation{\vecX}$ are
\begin{align*}
\GenMomentumFormg{\vecX}=& \insertion{\vecX} \StressEnergy &
\GenMomentum{\vecX}_\ia=&\StressEnergy_{\ia\ib}\vecX^{\ib} \\
\Deformation{\vecX}(\vecY,\vecZ)
=& (\CoDeriv_\vecY\formXg)(\vecZ)+(\CoDeriv_\vecZ\formXg)(\vecY)&
\Deformation{\vecX}_{\ia\ib}=& \CoDeriv_{\ia}\vecX_\ib + \CoDeriv_{\ib}\vecX_\ia ,
\end{align*}
which are related by Stokes' theorem
\begin{align*}
&&\int_{\partial\Omega} \GenMomentum{\vecX}_{\ia} \dThreeVec^{\ia}
=& \frac12 \int_\Omega \Deformation{\vecX}_{\ia\ib}\StressEnergy^{\ia\ib} \dFourGen .
\end{align*}
This is particularly useful for a Killing vector field, for which $\Deformation{\vecX}=0$. 

For any vector field $\vecX$, we will define the corresponding energy to be the hyper-surface integral of the generalized momentum
\begin{align*}
\GenEnergy{\vecX}[\MaxF](\sfc)
=\int_\sfc \GenMomentum{\vecX}_{\ia} \dThreeVec^{\ia} .
\end{align*}
This depends on the Maxwell field $\MaxF$ through the energy-momentum tensor. Frequently, we will be interested in $t=\text{const}$ hyper-surfaces, for which we define 
\begin{align*}
\GenEnergy{\vecX}[\MaxF](t)
=&\intTslice{t} \GenMomentum{\vecX}_{\ia} \dThreeVec^{\ia}
=\intTslice{t} \GenMomentum{\vecX}_{\ia} \vecTunit^\ia \dThreenTunit . 
\end{align*}
When the deformation tensor vanishes, by integrating over a space-time slab, one gets a conserved quantity: 
\begin{align*}
\GenEnergy{\vecX}[\MaxF](t_2)-\GenEnergy{\vecX}[\MaxF](t_1)
=0 .
\end{align*}
In applying this integration by parts, we require decay as $\rs\rightarrow\infty$, but merely smoothness as $\rs\rightarrow-\infty$, since in the maximal extension of the Schwarzschild manifold, for all values of $t$, the limit $r\rightarrow2M$, tends towards the same limiting sphere, the bifurcation sphere.

\subsection{Quantitative effect of trapping}
The energy-momentum tensor for the Maxwell field is
\begin{align}
\StressEnergy_{\ia\ib}=& \MaxF_{\ia\ic}\MaxF_\ib{}^{\ic} -\frac14 \g_{\ia\ib} \MaxF^{\ic\id}\MaxF_{\ic\id} \\
=& \frac12\left( \MaxF_{\ia\ic}\MaxF_\ib{}^{\ic} + (*\MaxF)_{\ia\ic}(*\MaxF)_\ib{}^{\ic} \right) . 
\label{eMaxEMS}
\end{align}
It satisfies \eqref{eEMSDivFree} and is trace-free
\begin{align*}
\g_{\ia\ib} \StressEnergy^{\ia\ib}=0 .
\end{align*}
Formally, one may assume that the Maxwell field is generated from a vector potential $\MaxA_\ia$ by $\MaxF=\Extd\MaxA$, $\MaxF_{\ia\ib}=\CoDeriv_\ia\MaxA_\ib-\CoDeriv_\ib\MaxA_\ia$ and take the Lagrangian to be $\Lagrangian= (1/2) \MaxF^{\ic\id}\MaxF_{\ic\id}=2(\CoDeriv_\ic\MaxA_\id)(\CoDeriv^\ic\MaxA^\id)$, in which case, the Lagrangian theory for the field $\MaxA$ gives the Maxwell equations \eqref{eMaxwellEquationDiv} as the Euler-Lagrange equations\footnote{The other equation, \eqref{eMaxwellEquationAlt}, holds because $\Extd^2=0$. } and \eqref{eMaxEMS} as the energy-momentum, which satisfies \eqref{eEMSDivFree}. Unfortunately, not all Maxwell fields can be represented in this way as an exterior derivative\footnote{The ``magnetically charged solution'', $F=q_{B} \sin(\theta) d\theta\wedge d\phi$ is not an exterior derivative. See appendix \ref{sExclusionOfNonRadiatable}. }. However, by direct computation from the Maxwell equation, it follows that the energy-momentum tensor in \eqref{eMaxEMS} satisfies \eqref{eEMSDivFree} so that Stokes' theorem can still be applied. 

The energy-momentum tensor is strictly positive when evaluated on time-like vectors. We will mainly be interested in time-like vectors with no angular components. Since any time-like vector with no angular component is a linear combination of $\veclunit$ and $\vecnunit$, to show the positivity of the stress-energy tensor, it is sufficient to compute the components in these null directions. These components are
\begin{align}
\StressEnergy(\veclunit,\veclunit)=& |\MyFp|^2 , \\
\StressEnergy(\veclunit,\vecnunit)=& |\MyF |^2 , \\
\StressEnergy(\vecnunit,\vecnunit)=& |\MyFm|^2 .
\label{eStressEnergyComponents}
\end{align}

The Schwarzschild manifold is static, so there is a conserved energy. The energy associated to the generator of $t$-translation, $\dt$, is strictly positive, 
\begin{align*}
\GenEnergy{\vecT}[\MaxF](t)
=&\frac12 \intTslice{t} \left(|\Elec|^2 +|\Mag|^2\right) \dThreeHF  \\
=&\frac14 \intTslice{t} \left(|\CKalpha|^2 +2|\CKrho|^2 +2|\CKsigma|^2 +|\CKalphab|^2\right) \dThreeHF  \\
=&\frac14 \intTslice{t} \left(|\MyFp|^2 +2|\MyF|^2 +|\MyFm|^2\right) \dThreeHF .
\end{align*}
The corresponding deformation tensor is 
\begin{align*}
\CoDeriv\formTg 
=&-\horifac^{-1}\formTg\otimes(-\frac{M}{r^2})\formRg +\horifac^{-1}\formRg\otimes(-\frac{M}{r^2})\formTg , \\
\Deformation{\vecT}
=& 0 .
\end{align*}
(We will need to compute deformation tensors later, but, in this case, we could simply have argued that the deformation tensor must vanish since $\dt$ generates a symmetry of the metric.) From the vanishing of the deformation tensor, we have a conservation law
\begin{align*}
\GenEnergy{\vecT}[\MaxF](t)
=&\GenEnergy{\vecT}[\MaxF](0) .
\end{align*}
This immediately gives an upper bound on the average value of the components of the Maxwell tensor in any region bounded away from the event horizon. 

By applying Lie derivatives, we can get additional conservation laws. If $\vecX$ generates a symmetry and $\MaxF$ solves the Maxwell equations, then $\Lie_\vecX \MaxF$ will also be a solution of the Maxwell equations. For each set of symmetries and integer $k$, we have the conserved quantities
\begin{align*}
\GenEnergy{\vecT}[\Lie_\AngAlgGen^k\MaxF](t)=&\GenEnergy{\vecT}[\Lie_\AngAlgGen^k\MaxF](0) ,\\
\GenEnergy{\vecT}[\Lie_\SymGen^k\MaxF](t)=&\GenEnergy{\vecT}[\Lie_\SymGen^k\MaxF](0) .
\end{align*}

We now improve these estimates and reveal the effect of trapping, by considering the conformal energy. Following earlier work \cite{CK,BlueSofferLongPaper,BlueSterbenz,DafermosRodnianski}, we let
\begin{align*}
\vecConf=& (t^2+\rs^2)\dt +2t\rs\dr \\
=&\frac12\left(\uout^2\veclcoord +\uin^2\vecncoord \right).
\end{align*}
We will call this the conformal vector field, but it is also one of the vector fields known as the Morawetz vector field. It is an analogue of a vector field used in $\Reals^{1+n}$ to prove decay for the wave equation, the Maxwell equation, and the spin $2$ field. The analogue in $\Reals^{1+n}$ generates a positive quantity, so it is not surprising that the same holds on the Schwarzschild manifold. We define the conformal energy to be
\begin{align}
\GenEnergy{\vecConf}[\MaxF](t)
=& \intTslice{t} \GenMomentum{\vecConf}_\ia \dThreeVec^{\ia} \nonumber\\
=& \intTslice{t} \left((1/2)(t^2+\rs^2)(|\Elec|^2+|\Mag|^2) 
+2t\rs(\Elec_\vecHunit\Mag_\vecPhunit-\Elec_\vecPhunit\Mag_\vecHunit)\right) \dThreeHF\nonumber \\
=& (1/4)\intTslice{t} \left(
(t+\rs)^2|\CKalpha|^2 +2(t^2+\rs^2)(|\CKrho|^2+|\CKsigma|^2) +(t-\rs)^2|\CKalphab|^2\right) \dThreeHF \nonumber \\
=& (1/4)\intTslice{t} \left( \uout^2|\MyFp|^2 + (\uout^2+\uin^2)|\MyF|^2 +\uin^2|\MyFm|^2\right) \dThreeHF . \nonumber
\end{align}
In the null decomposition or spinor representation, all the terms in the integrand are non-negative, and, inside the light-cone $|\rs|<(1-\epsilon)t$, the coefficients on the Maxwell field components grow like $t^2$. Thus, once we show that the conformal energy is bounded, there will be decay for the localized field components. 

The following lemma gives an almost conservation law for the conformal energy. It states that, to bound the conformal energy, it is sufficient to prove sufficiently strong decay in a particularly region bounded away from the event horizon. There are two important observations to make from this lemma and its proof: (i) an estimate for the two field components $\Elec_\vecRunit$ and $\Mag_\vecRunit$ will control all the field components through the conformal charge, and (ii) it is sufficient to control these field components only in a region near the photon sphere $r=3M$. 

\begin{lemma}[Trapping lemma]
There is a positive function $\chiTrap$ supported in a bounded range of $\rs$ values such that if $\MaxF$ is a solution to the Maxwell equations \eqref{eMaxwellEquationDiv}-\eqref{eMaxwellEquationAlt}, then
\begin{align}
\GenEnergy{\vecConf}[\MaxF](t_2)-\GenEnergy{\vecConf}[\MaxF](t_1)
\leq& \intTslab{t_1}{t_2} t\chiTrap  \left(\Elec_\vecRunit^2 +\Mag_\vecRunit^2\right) \dFourEx \nonumber\\
\leq& \intTslab{t_1}{t_2} t\chiTrap  |\MyF|^2 \dFourEx .
\label{eTrappingEstimate}
\end{align}
\begin{proof}
The deformation tensor for $\vecConf$ is given by 
\begin{align*}
\CoDeriv\formRg
=&r^{-1}\horifac\g - r^{-1}(1-3M/r)(-\formTunitg\otimes\formTunitg+\formRunitg\otimes\formRunitg ) ,\\
\CoDeriv\formConfg
=&(t^2+\rs^2)\CoDeriv\formTg + 2t\rs\CoDeriv\formRg \nonumber\\
&-\horifac^{-1}\formTg\otimes 2t\formTg +\horifac^{-1}\formRg\otimes 2\rs\formTg\nonumber\\
&-\horifac^{-1}\formTg\otimes 2\rs\formRg +\horifac^{-1}\formRg\otimes 2t\formRg ,\\
\Deformation{\vecConf}
=&2t\rs \Deformation{\vecR}
+4t (-\formTunitg\otimes\formTunitg +\formRunitg\otimes\formRunitg)\nonumber\\
=&4t\frac{\rs}{r}\horifac\g + 4t\left(1-\frac{\rs}{r}(1-3M/r)\right)(-\formTunitg\otimes\formTunitg+\formRunitg\otimes\formRunitg ) .
\end{align*}
Because the Maxwell energy-momentum tensor is trace-free, the contraction of the first term against $\StressEnergy$ is zero at each point. The importance of $r=3M$, where the orbiting geodesics are located, is immediate from the second term. The contraction against the energy-momentum tensor is
\begin{align*}
\Deformation{\vecConf}_{\ia\ib}\StressEnergy^{\ia\ib}
=& 4t\left(\frac{\rs}{r}(1-3M/r)-1\right)(\StressEnergy_{\vecTunit\vecTunit}-\StressEnergy_{\vecRunit\vecRunit}) , \\
=& 4t\left(\frac{\rs}{r}(1-3M/r)-1\right) \StressEnergy(\veclunit,\vecnunit) .
\end{align*}
From this, we have the following almost-conservation law
\begin{align}
\GenEnergy{\vecConf}[\MaxF](t_2)-\GenEnergy{\vecConf}[\MaxF](t_1)
=& 2\intTslab{t_1}{t_2} t\left(1-\frac{\rs}{r}(1-3M/r)\right) \left(\Elec_\vecRunit^2 +\Mag_\vecRunit^2\right) \dFourEx \nonumber \\
=& 2\intTslab{t_1}{t_2} t\left(1-\frac{\rs}{r}(1-3M/r)\right) |\MyF|^2 \dFourEx .
\label{eMaxConfEstimate}
\end{align}
We refer to $1-\frac{\rs}{r}(1-3M/r)$ as the trapping term. 

For $r\rightarrow 2M$, $\rs\rightarrow-\infty$ and $1-3M/r\rightarrow -1/2$, so $1-(1-3M/r)\rs/r$ is negative. The explicit expression for $\rs$ in terms of $r$ is 
\begin{align*}
\rs=& r +2M\log\left(\frac{r-2M}{2M}\right) - 3M +2M\log2 .
\end{align*}
Because of the logarithmic term, as $r\rightarrow\infty$, $(1-\frac{\rs}{r}(1-3M/r))=(r-\rs)/r +O(1/r) < -2M \log(r) +O(1/r)$ which is negative for sufficiently large $r$. Since the trapping term has negative limit at $\pm\infty$, it is positive only in a compact interval. 

We now introduce a smooth, compactly supported function $\chiTrap$ which dominates the trapping term. This function is chosen to satisfy $4\left((1-3M/r)\rs/r-1 \right)<\chiTrap$. This gives the desired result. 
\end{proof}
\end{lemma}

\subsection{Spin reduction}
In this section, we obtain a decay result for the zero-weight component. From the previous section, we know this is enough to control energies involving all the components. It is known that the evolution of the zero-weight component can be determined from a wave equation without referring to the other components. Thus, we can reduce the problem from the Maxwell equations to a wave equation. Since physicists refer to wave equations as spin $0$ equations and the system of the Maxwell equations as a spin $1$ system, we use ``spin reduction'' to refer to this reduction. 

The Maxwell equations can be written as a fairly simple system in terms of the null coordinate bases and the corresponding components. This is a result due to Price \cite{Price}, although, he uses a null tetrad, which makes the corresponding expressions look significantly different. By direct computation and application of the Maxwell equations, 
\begin{align}
\vecncoord\CoordWtp =& \vecmcoord \CoordWt \horifac r^{-2} ,\label{ePricei}\\
\veclcoord\CoordWt  =& \vecmbcoord\CoordWtp +\cot\theta \CoordWtp ,\label{ePriceii}\\
\vecncoord\CoordWt  =& -\vecmcoord\CoordWtm -\cot\theta \CoordWtm ,\label{ePriceiii}\\
\veclcoord\CoordWtm =& -\vecmbcoord \CoordWt \horifac r^{-2} . \label{ePriceiv}
\end{align}
We refer to these as the ``Price equations''. 

The cotangent terms appear to be singular; however, if $\CoordWtp$ is treated as spherical $1$-forms, then the combination of the angular derivative and the cotangent term can be written simply as 
\begin{align}
\vecmbcoord\CoordWtp +\cot\theta \CoordWtp
=& \horifac^{-1/2}\left(\text{div}\CKalpha + i\text{curl}\CKalpha \right), \\
\vecmcoord\CoordWtm +\cot\theta \CoordWtm
=& \horifac^{-1/2}\left(\text{div}\CKalphab - i\text{curl}\CKalphab \right), 
\label{eMPlusCotIsDivCurl}
\end{align}
where $\text{div}$ and $\text{curl}$ are the spherical divergence and curl. If we'd defined a coordinate based null decomposition $A(\vecea)=F(\veclcoord,\vecea)=\horifac^{-1/2}\CKalpha(\vecea)$, then we'd have exactly $\vecmbcoord\CoordWtp +\cot\theta \CoordWtp = \text{div} A +i\text{curl} A$, and similarly for the other components. One important consequence of this is that the right-hand sides of \eqref{ePriceii} and \eqref{ePriceiii} are controlled by
\begin{align*}
|\vecmbcoord\CoordWtp +\cot\theta \CoordWtp|
+|\vecmcoord\CoordWtm +\cot\theta \CoordWtm|
\leq& \horifac^{-1/2} r \TensorDNorm{\MaxF}{\NiceNormalisedBasisSet}{1}{\AngAlgGen} .
\end{align*}

Another important consequence of the Price equations is that the zero weight term satisfies a wave equation. From \eqref{ePriceiii} and \eqref{ePriceiv}, 
\begin{align*}
\veclcoord\vecncoord \CoordWt =& (\vecmcoord+\cot\theta)\vecmbcoord\CoordWt \horifac r^{-2} ,\\
-\dtt \CoordWt =& -\drr\CoordWt +r^{-2}\horifac(-\slap)\CoordWt .
\end{align*}
If there were an additional $(2M/r^3)\horifac \CoordWt$ term on the right, then $\CoordWt$ would be a solution to the wave equation on the Schwarzschild manifold, $\CoDeriv^\ia\CoDeriv_\ia (r^{-1}\CoordWt)$. Even in the absence of this term, the previous analysis of wave equations is sufficiently general to apply to a wave equation of this form \cite{BlueSterbenz}. In fact, the wave equation under consideration is simpler than the true wave equation\footnote{If $\CoDeriv^\ia\CoDeriv_\ia(r^{-1}\SFn)=0$, then $\SFn$ satisfies $-\dtt\SFn= -\drr\SFn +V\SFn +V_L(-\slap)\SFn$ with $V=2M r^{-3} \horifac$, thus, the equation governing $\CoordWt$ is closer to $\CoDeriv^\ia\CoDeriv_\ia(r^{-1}\SFn)=0$ than to $\CoDeriv^\ia\CoDeriv_\ia \SFn=0$. } $\CoDeriv^\ia\CoDeriv_\ia(r^{-1}\SFn)=0$, and, in appendix \ref{sOneDWaveAnalysis}, we provide a stream-lined version of the method from \cite{BlueSterbenz}.

For solutions to a wave equation, there are estimates on the weighted space-time integral we need to control for the conformal estimate. If $\SFn$ is a solution to
\begin{align}
-\dtt\SFn=& -\drr\SFn + V_L(-\slap)\SFn
\label{ePriceWaveEqn}
\end{align}
with
\begin{align*}
V_L=r^{-2}\horifac ,
\end{align*}
then the energy and conformal charge are defined in terms of their densities by
\begin{align*}
\PriceEnergyDensity=& |\dt\SFn|^2 + |\dr\SFn|^2 + V_L |\dAng\SFn|^2 ,\\
\PriceConfDensity=& \frac14 |(t+r)(\dt+\dr)\SFn|^2 
+\frac14|(t-r)(\dt-\dr)\SFn|^2
+\frac12 (t^2+\rs^2)V_L|\dAng\SFn|^2 +\PriceEnergyDensity ,\\
\PriceEnergy[u](t)=& \frac12 \intTslice{t} \PriceEnergyDensity \dThreePOneD ,\\
\PriceConf[u](t)=& \frac12 \intTslice{t} \PriceConfDensity \dThreePOneD .
\end{align*}
The energy and conformal energy are generated from $\vecT$ and $\vecConf$. (In fact, the conformal energy is generated by $\vecConf+\vecT$ to provide better control at time $t$ near $0$.) Since $\vecT$ generates a symmetry, the energy is conserved. As with the conformal energy for the Maxwell field, the conformal energy is not conserved, and the energy density near the photon sphere $r=3M$ must be controlled. The important results for this discussion are that, at any time $t>0$, $k\in\Naturals$, and any compactly supported function $\chi$, 
\begin{align}
\PriceEnergy[\SFn](t)=&\PriceEnergy[\SFn](0) ,\nonumber\\
\PriceEnergy[\dAng^k\SFn](t)=&\PriceEnergy[\dAng^k\SFn](0) ,\nonumber\\
\PriceConf[\SFn(t)]\leq& \PriceConf[\SFn(0)] + C \PriceEnergy[\slap^2\SFn(0)] ,\nonumber\\
\intTslabInfinity{0} \frac{|\SFn|^2}{(1+\rs^2)^2} \dFourPOneD \leq& \PriceEnergy[\SFn](0) ,\label{eWaveLocalDecay} \\
\intTslabInfinity{0} t \chi |\dAng\SFn|^2 \dFourPOneD \leq& C \PriceConf[\SFn](0) +\PriceEnergy[\slap^2\SFn](0). \nonumber
\end{align}
In appendix \ref{sExclusionOfNonRadiatable}, we exclude spherically symmetric harmonics, so, from dropping the angular derivatives on the left-hand side of the previous estimate, 
\begin{align}
\intTslabInfinity{0} t \chiTrap|\SFn|^2 \dFourEx \leq&
C\left( \PriceConf[\SFn](0) + C \PriceEnergy[\slap^2\SFn](0) \right).\label{eWaveEqnTrappingEstimate}
\end{align}

We apply this with $\SFn=\CoordWt$ in the following lemma. 

\begin{lemma}
\label{lSurfaceEnergyBound}
If $\MaxF$ is a solution to the Maxwell equations \eqref{eMaxwellEquationDiv}-\eqref{eMaxwellEquationAlt} and $n\geq0$, then
\begin{align}
\GenEnergy{\vecConf}[\Lie_\SymGen^n \MaxF](t)
\leq& C \left(\normKangN{n+1}+\normTangN{n+5}\right) , \label{eConformalEnergyBound}\\
\intTslabInfinity{0} t\chiTrap |\MyF|^2 \dFourEx 
\leq& C\left( \sum_{k=0}^1\GenEnergy{\vecConf}[\Lie_\AngAlgGen^k \MaxF](0) +\sum_{k=0}^{5}\GenEnergy{\vecT}[\Lie_\AngAlgGen^k \MaxF](0) \right) .
\label{eGlobalConfEstimate} 
\end{align}
Furthermore, if the normal to $\sfc$ has uniformly bounded below $\veclunit$ and $\vecnunit$ components, then
\begin{align*}
\int_\sfc \left( |\MyFp|^2 +2|\MyF|^2 +|\MyFm|^2\right) \dThreeHF
\leq& C \max_\sfc(\uin^{-2},\uout^{-2}) \left(\normKangN{1}+\normTangN{5}\right) .
\end{align*}
\begin{proof}
Taking $\CoordWt=\SFn$, the energy associated to the Maxwell field $\MaxF$ and that of the scalar wave $\SFn$ are closely related. From the Price equations \eqref{ePricei}-\eqref{ePriceiv} and the geometric interpretation of the $\vecmcoord+\cot\theta$ terms in \eqref{eMPlusCotIsDivCurl}, 
\begin{align}
\PriceEnergy[\SFn](t)=& \GenEnergy{\vecT}[\Lie_\AngAlgGen \MaxF](t) , \label{ePrivevsTrueEnergy}\\
\PriceEnergy[\slap^2\SFn](t)=& \GenEnergy{\vecT}[\Lie_\AngAlgGen^5 \MaxF](t) , \label{ePrivevsTrueEnergyDeriv}\\
\PriceConf[\SFn](t)=& \GenEnergy{\vecConf}[\Lie_\AngAlgGen \MaxF](t) \label{ePrivevsTrueConfEnergy}. 
\end{align}
Estimate \eqref{eWaveEqnTrappingEstimate} can be written as
\begin{align*}
\intTslabInfinity{0} t\chiTrap |\MyF|^2 \dThreeHF
\leq& C( \GenEnergy{\vecConf}[\Lie_\AngAlgGen \MaxF](0) +\sum_{k=0}^{5}\GenEnergy{\vecT}[\Lie_\AngAlgGen^k \MaxF](0) ) .
\end{align*}
From this estimate and the trapping estimate \eqref{eTrappingEstimate}, the estimate \eqref{eConformalEnergyBound} follows. If one of the derivatives in \eqref{eConformalEnergyBound} is in the angular direction instead of the time direction, then it would not be necessary to drop the angular derivative in \eqref{eWaveEqnTrappingEstimate}, and only $n$ and $n+4$ derivatives would be needed on the $\vecConf$ and $\vecT$ energies respectively. 

Since the integral of the trapping term has been controlled in the entire of the exterior of the Schwarzschild manifold, we have a uniform bound on the integral of $\GenMomentum{\vecConf}$ on any hyper-surface. If the hyper-surface has a normal with uniformly  bounded below $\veclunit$ and $\vecnunit$ components, then the integral will be bounded below by $C(\uout^2|\MyFp|^2 +(\uin^2+\uout^2)|\MyF|^2 +\uin^2|\MyFm|^2)$. This provides the final estimate. 
\end{proof}
\end{lemma}

We remark that from equations \eqref{ePrivevsTrueEnergy}-\eqref{ePrivevsTrueConfEnergy}, we could have bounded the energies $\GenEnergy{\vecT}[\Lie_\AngAlgGen\MaxF](t)$ and $\GenEnergy{\vecConf}[\Lie_\AngAlgGen \MaxF](t)$ by immediately appealing to results for the wave equation. However, we would still need to present the energy and conformal energy for the spin $0$ wave equation and for the spin $1$ Maxwell equation, relate them, and present the Price equations. While it would have been possible for us to omit the Lagrangian theory and the trapping lemma, this would have removed the motivation for considering the energy and conformal energy.

%
%
\section{Pointwise decay in stationary regions}
\label{sStationaryDecay}
In this section, our goal is to prove $L^\infty$ decay in regions where $2M<r_1<r<r_2$. We refer to these as stationary regions since the range of the radial coordinate does not change in $t$. Restricting attention to a stationary region, the integrand in the conformal energy behaves like $t^2$ times the Maxwell field components squared. Since the conformal energy is bounded, the field components decay in $L^2_{\text{loc}}$ like $t^{-1}$. 

Control on radial derivatives is the main thing that we need to improve this from decay in mean to pointwise decay. Sobolev estimates can be used to convert $L^2_{\text{loc}}$ decay for derivatives into $L^\infty_{\text{loc}}$ decay. For this, we need decay on the spatial derivatives of the Maxwell field. From spherical symmetry, the Lie derivative of the Maxwell field in the direction of an angular derivative, $\Lie_{\RotElementi}\MaxF$, also satisfies the Maxwell equations and has the same type of decay in mean as $\MaxF$. Since $\vecR$ does not generate a symmetry, the Lie derivative in that direction will not solve the Maxwell equations. 

To control the radial derivatives, we use the structure of the Maxwell equations. Using the staticity of the Schwarzschild manifold, we can control $t$ derivatives, $\Lie_{\vecT}\MaxF$, in $L^2_{\text{loc}}$. In a fixed, compact range of $r$ values, the covariant derivatives of the coordinate basis are controlled by finite multiples of the coordinate bases again. We are working in $L^2$ where we already control all the components. Thus, we control the difference between components of the covariant derivative in a direction and the covariant derivative of the components of the Maxwell tensor (ie, $\CoDeriv_\ia \MaxF_{\ib\ic}\vecX^\ib\vecY^\ic\sim \CoDeriv_\ia(\MaxF_{\ib\ic}\vecX^\ib\vecY^\ic)$). 

The notation in subsection \ref{ssLieDerivatives} can be used to define a ``big-$O$'' notation to estimate the difference between two functions depending on position and a tensor field. We say a function of position and a $(0,2)$-tensor field is equal to another such function up to norm terms and in an interval, if, on any bounded interval of $\rs$ values, there is a constant such that, for any $(0,2)$ tensor, the difference between the two functions is bounded by a constant times the norm of the tensor
\begin{align*}
f=&h +\ErrorTermsOf{\tensorA} &\Longleftrightarrow&
&|f(t,\rs,\theta,\phi,\tensorA) - h(t,\rs,\theta,\phi,\tensorA)| \leq& C\TensorNorm{\tensorA}{\NiceBasisSet} .
\end{align*}
Similarly, for two collections of such functions, we say
\begin{align*}
\{ f \} =& \{ h \} +\ErrorTermsOf{\tensorA} 
\end{align*}
if for each $f$ there is an $h$ such that $f=h+\ErrorTerms{\tensorA}$ and vice versa. We say
\begin{align*}
\{ f \} \lesssim& \{ h \} +\ErrorTermsOf{\tensorA} 
\end{align*}
if, for any bounded interval in $\rs$, there is a constant $C$ such that each $|f|$ is bounded by $C$ times the sum of the absolute values of the $h$'s plus $\ErrorTermsOf{\tensorA}$ terms. We make similar definitions involving $\ErrorTermskOf{\tensorA}{k}$. 

The big-$O$ notation used here is local to compact intervals, which allows us to ignore the difference between normalized and unnormalized vector fields, 
\begin{align*}
\TensorNorm{\tensorA}{\NiceNormalisedBasisSet}=&\ErrorTermsOf{\tensorA} ,\\
\TensorNorm{\tensorA}{\NiceBasisSet}=&\BigO{\TensorNorm{\tensorA}{\NiceNormalisedBasisSet}} .
\end{align*}

This notation allows us to prove theorem \ref{tDecayInStationaryRegions}. As outlined in the beginning of this section, our strategy in the proof is to use the Maxwell equations to trade derivatives along the generators of symmetries for spatial derivatives and then to apply the Sobolev estimate. In doing this, we use the big-$O$ notation to estimate error terms generated by converting between Lie and covariant derivatives. This allows us to improve our decay estimates from decay in mean to pointwise decay. Here, we use the null decomposition and explicitly state the norms. Clearly the same result holds for $|\Elec|+|\Mag|$ or $|\CKalpha|+|\CKrho|+|\CKsigma|+|\CKalphab|$. 

\begin{theorem}
\label{tDecayInStationaryRegionsInText}
Let $2M<r_1<r_2<\infty$. There is a constant $C_{(r_1,r_2)}$ such that if $\MaxF$ is a solution of the Maxwell equations \eqref{eMaxwellEquationDiv}-\eqref{eMaxwellEquationAlt}, then for all $t\in\Reals$, $r\in[r_1,r_2]$, and $(\theta,\phi)\in S^2$,
\begin{align*}
|\MyFp|+|\MyF|+|\MyFm| \leq& C_{(r_1,r_2)} t^{-1} \left(\normKangN{4} +\normTangN{8}\right)^{1/2} .
\end{align*}
\begin{proof}
Using the big-$O$ notation, we can control the difference between the derivative of a component of the Maxwell field and the corresponding component of the Lie derivative. Since the Lie derivative of any vector field with respect to any other is a linear combination of the coordinate vector fields with smooth coefficients, 
\begin{align*}
\Lie_\NiceBasisSet (\MaxF(\NiceBasisSet,\NiceBasisSet))
=& (\Lie_\NiceBasisSet\MaxF)(\NiceBasisSet,\NiceBasisSet) + \ErrorTerms .
\end{align*}
This process can be iterated, so that 
\begin{align*}
\Lie_{\SymGen}^k(\MaxF(\NiceBasisSet,\NiceBasisSet)) =& \ErrorTermsk{k} .
\end{align*}
Similarly for covariant derivatives, 
\begin{align*}
\NiceBasisSet (\MaxF(\NiceBasisSet,\NiceBasisSet))
=& (\CoDeriv_\NiceBasisSet\MaxF)(\NiceBasisSet,\NiceBasisSet) + \ErrorTerms .
\end{align*}
We note that if we had applied two symmetry-generating derivatives before making the estimate we would have
\begin{align*}
\Lie_\SymGen\Lie_\SymGen( \NiceBasisSet (\MaxF(\NiceBasisSet,\NiceBasisSet)) )
=& \Lie_\SymGen\Lie_\SymGen ((\CoDeriv_\NiceBasisSet\MaxF)(\NiceBasisSet,\NiceBasisSet) )
+ \ErrorTermsk{2} ,
\end{align*}
and similarly with the Lie derivative replacing the covariant derivative in $\NiceBasisSet$. 

To control the radial derivative of components which have no $\vecR$ arguments, we use \eqref{eMaxwellEquationAlt}, 
\begin{align*}
\Lie_\vecR \Lie_\SymGen\Lie_\SymGen (\MaxF(\SymGen,\SymGen))
=& \Lie_\SymGen\Lie_\SymGen \Lie_\vecR(\MaxF(\SymGen,\SymGen)) \\
=& \Lie_\SymGen\Lie_\SymGen \CoDeriv_\vecR(\MaxF(\SymGen,\SymGen)) \\
=& \Lie_\SymGen\Lie_\SymGen ((\CoDeriv_\vecR\MaxF)(\SymGen,\SymGen))
+ \ErrorTermsk{2} \\
=& \Lie_\SymGen\Lie_\SymGen ((\CoDeriv_\SymGen\MaxF)(\vecR,\SymGen))
+ \ErrorTermsk{2} \\
=& \Lie_\SymGen\Lie_\SymGen \CoDeriv_\SymGen(\MaxF(\vecR,\SymGen))
+ \ErrorTermsk{2} \\
=& \Lie_\SymGen\Lie_\SymGen \Lie_\SymGen(\MaxF(\vecR,\SymGen))
+ \ErrorTermsk{2} \\
=& \ErrorTermsk{3} .
\end{align*}
Similarly, to gain control of component with one radial argument, we apply \eqref{eMaxwellEquationDiv}, 
\begin{align*}
\Lie_\vecR \Lie_\SymGen\Lie_\SymGen (\MaxF(\vecR,\SymGen))
\lesssim& \Lie_\SymGen\Lie_\SymGen ((\CoDeriv_\vecRunit\MaxF)(\vecRunit,\SymGen))
+ \ErrorTermsk{2} \\
\lesssim& \Lie_\SymGen\Lie_\SymGen ((\CoDeriv_\SymGen\MaxF)(\SymGen,\SymGen))
+ \ErrorTermsk{2} \\
=& \ErrorTermsk{3} .
\end{align*}
Since $\MaxF$ is antisymmetric, there is no need to control components with two $\vecR$ arguments. 

Control of triple derivative terms, of the form $\Lie_\NiceBasisSet\Lie_\SymGen\Lie_\SymGen (\MaxF(\NiceBasisSet,\NiceBasisSet))$, is sufficient to prove $L^\infty$ decay. From the boundedness of the conformal charge, for any interval $[r_1,r_2]$ in the exterior, there is a constant $C$, such that (with ${\rs}_1$ and ${\rs}_2$ the values of $\rs$ corresponding to $r=r_1$ and $r=r_2$) 
\begin{align*}
\intTBox{t}{{\rs}_1}{{\rs}_2} t^2 \TensorNorm{\MaxF}{\NiceBasisSet}^2 \dThreePOneD
\leq& C\GenEnergy{\vecConf}[\MaxF] .
\end{align*}
A local, inhomogeneous, $1$-dimensional Sobolev estimate gives
\begin{align*}
\int_{\{t\}\times\{\rs\}\times S^2} |\Lie_\SymGen\Lie_\SymGen (\MaxF(\NiceBasisSet,\NiceBasisSet))|^2 \dTwo 
\leq& C \intTBox{t}{{\rs}_1}{{\rs}_2} |\Lie_\vecR \Lie_\SymGen\Lie_\SymGen (\MaxF(\NiceBasisSet,\NiceBasisSet))|^2 +|\Lie_\SymGen\Lie_\SymGen (\MaxF(\NiceBasisSet,\NiceBasisSet))|^2 \dThreePOneD \\
\leq& C \intTBox{t}{{\rs}_1}{{\rs}_2}  \TensorDSymNorm{\MaxF}{3}^2 \dThreePOneD \\
\leq& Ct^{-2} \sum_{k=0}^{3}\GenEnergy{\vecConf}[\Lie_\SymGen^3\MaxF] .
\end{align*}
Now applying a spherical Sobolev estimate, we have
\begin{align*}
|\MaxF(\NiceBasisSet,\NiceBasisSet)(t,\rs,\theta,\phi) |
\leq& Ct^{-1} \left(\sum_{k=0}^{3}\GenEnergy{\vecConf}[\Lie_\SymGen^k \MaxF](t)\right)^\frac12 .
\end{align*}

By lemma \ref{lSurfaceEnergyBound}, the conformal energy at any time is bounded by the initial data (with extra derivatives). This gives, in any stationary region away from the event horizon, that the components decay like $t^{-1}$. 
\end{proof}
\end{theorem}

%
%
\section{Decay outside stationary regions}
\label{sMovingDecay}
In this section, we prove decay for the field components outside of stationary regions. In Minkowski space $\Reals^{1+3}$, it is typical to obtain decay estimates in the regions $|\vec{x}|<(1-\epsilon)t$ and $|\vec{x}|>(1-\epsilon)t$. Because boosts are not symmetries of the Schwarzschild solution, decay in a stationary region is different from decay along outgoing curves $\rs\sim(1-\epsilon) t$. Similarly, since there is no reflection symmetry $\rs\rightarrow-\rs$, the decay rates in the regions $\rs>0$ and $\rs<0$ are different. Thus, we obtain decay in stationary regions, outgoing regions, and ingoing regions. 

Most of the decay estimates in this section are proven by considering the energy on ingoing or outgoing null hyper-surfaces. We will use $\NullInSfc{\uout}$ and $\NullOutSfc{\uin}$ to refer to ingoing and outgoing null hyper-surfaces on which $\uout$ and $\uin$ are constant respectively. We will restrict these to the future $t\geq0$. To make estimates on $\NullOutSfc{\uin}$, we can introduce a parameter $t_1$ and an approximate surface which extends along the hyper-surface $t=0$ from the bifurcation sphere ($\rs\rightarrow-\infty)$ to the intersection of $t=0$ with $\NullOutSfc{\uin}$, extends along $\NullOutSfc{\uin}$ until $t=t_1$, and finally continues onto space-like infinity along $t=t_1$. Since the deformation tensor for $\vecT$ is zero, the surface integral of the generalized momentum $\GenMomentum{\vecT}$ along this surface will be the same as the integral along $t=0$. Similarly, since estimate \eqref{eGlobalConfEstimate} says that the integral over the entire exterior region of the positive part of the $\vecConf$ deformation tensor is bounded, the surface integral of the generalized momentum $\GenMomentum{\vecConf}$ over the approximating surface is bounded by the initial data. Dropping the positive contribution from integrating along $t=0$ and $t=t_1$ and taking the limit as $t_1\rightarrow\infty$, 
\begin{align*}
\intOutSfc \GenMomentum{\vecConf}_\ia \dThreeVec^\ia \leq& C\left(\normKangN{1}+\normTangN{5}\right) . 
\end{align*}
A similar argument can be made for $\NullInSfc{\uout}$. Since $\NullOutSfc{\uin}$ and $\NullInSfc{\uout}$ are null surfaces, we cannot apply the last part of lemma \ref{lSurfaceEnergyBound}.

These integrals can be expanded in terms of the Maxwell field components. 
\begin{align*}
\intInSfc \left(\uout^2|\MyF|^2 +\uin^2|\MyFm|^2\right) \horifac r^2\dThreeUInPOD 
=&  \intInSfc \GenMomentum{\vecConf}_\ia \dThreeVec^\ia,\\
\intOutSfc \left(\uin^2|\MyF|^2 +\uout^2|\MyFp|^2\right) \horifac r^2\dThreeUOutPOD 
=& \intOutSfc \GenMomentum{\vecConf}_\ia \dThreeVec^\ia . 
\end{align*}
To obtain estimates for derivatives tangential to this surface, we can convert to the coordinate based components and apply angular derivatives and the Price equations \eqref{ePricei}-\eqref{ePriceiv} to get 
\begin{align}
\intInSfc \left(\uout^2|\vecncoord\CoordWtp|^2\horifac^{-1} r^{2} +\uin^2|\vecncoord\CoordWt|^2\right) \dThreeUInPOD 
\leq& C(\normKangN{4}+\normTangN{8}) ,\label{eLabelForAngledConformalIn}\\
\intOutSfc \left(\uin^2|\veclcoord\CoordWtm|^2\horifac^{-1} r^{2} +\uout^2|\veclcoord\CoordWt|^2\right) \dThreeUOutPOD 
\leq& C(\normKangN{4}+\normTangN{8}) .
\label{eLabelForAngledConformalOut}
\end{align}

We now prove decay in outgoing regions. 

\begin{lemma}[Decay for $\rs>1$]
\label{lFarDecay}
There is a constant $C$ such that if $\MaxF$ is a solution of the Maxwell equations \eqref{eMaxwellEquationDiv}-\eqref{eMaxwellEquationAlt}, then for all $t\geq0$, $\rs>1$, $(\theta,\phi)\in S^2$, 
\begin{align*}
|\MyFp|\leq& C r^{-3/2} |\uout|^{-1} \left( \normKangN{4} +\normTangN{8} +\SupInitData \right)^{1/2} , \\
|\MyF|\leq& C r^{-2} \left(\frac{\uout-|\uin|}{\uout(1+|\uin|)}\right)^{1/2} \left( \normKangN{4} +\normTangN{8} +\SupInitData \right)^{1/2}, \\
|\MyFm|\leq& C r^{-1} (1+|\uin|)^{-1} \left( \normKangN{4} +\normTangN{8} +\SupInitData \right)^{1/2} .
\end{align*}
For $t<\rs$, 
\begin{align*}
|\MyFm|\leq& C r^{-1} (1+|\uin|)^{-3/2} \left( \normKangN{4} +\normTangN{8} +\SupInitData \right)^{1/2} .
\end{align*}
\begin{proof}
At any point in the far region, $\rs>1$, we will integrate along a radial, null ray to prove decay. The bounds on the conformal charge give decay for integrals along the null rays. The final end point will either be at $t=0$, where we already have decay, or at $\rs=0$, where we have decay by assumption. In this way, each component of the Maxwell field will be written as the sum of two terms, both of which decay. The typical null rays which we use are illustrated in figure \ref{FigOuterNullRays}. 

There are a number of simplifications in the outgoing region. We can ignore factors of $\horifac$, since the ratio between $1$ and $\horifac$ is bounded above and bounded below by a strictly positive number. Since $\rs>1$, we can ignore ratios of $\rs/r$. There is the ordering $\uout\geq\rs\geq\uin$. On outgoing null rays, on which $\uin$ is constant, the change in $\uout$ is twice the change in $\rs$, and similarly, on ingoing radial, null rays the change in $\uin$ is twice the change in $\rs$. 

\begin{figure}
\begin{center}
\input{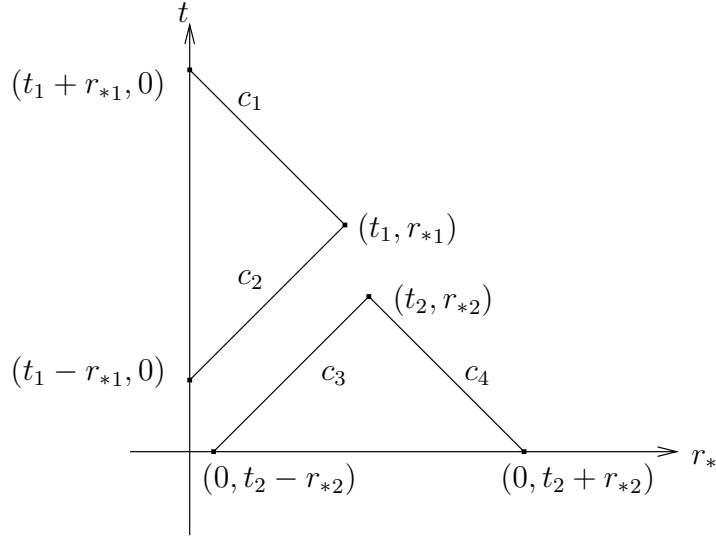}
\caption{Null rays in the outer region, $\rs>0$. The angular variables have been suppressed. The null rays go from a point either to the initial hypersurface $t=0$ or to the stationary region $\rs=0$. } 
\label{FigOuterNullRays}
\end{center}
\end{figure}

The simplest application of our method is for the zero-weight component. First, we prove an estimate inside the light-cone, for $t>\rs$. We use a radial, null geodesic from $(t,\rs,\theta,\phi)$ to $(t+\rs,0,\theta,\phi)$ parametrized by $\uin$, 
\begin{align*}
|\CoordWt(t,\rs,\theta,\phi)|
\leq& \int_{c_1} |\vecncoord\CoordWt| d\uin +|\CoordWt(t+\rs,0,\theta,\phi)|\\
\leq& \left(\int_{c_1} \uin^{-2} d\uin\right)^\frac12 
\left(\int_{c_1} \uin^2 |\vecncoord\CoordWt|^2 d\uin\right)^\frac12 
+ |\CoordWt(t+\rs,0,\theta,\phi)|\\
\leq& \left(\frac1{\uin}-\frac1{\uout}\right)^{1/2}\left(\int_{c_1} \uin^2 |\vecncoord\CoordWt|^2 d\uin\right)^\frac12 + |\CoordWt(t+\rs,0,\theta,\phi)|  .
\end{align*}
The end point decays at a rate of $t'^{-1}$ evaluated at $t'=t+\rs$. We now integrate over the angular variables too (and apply Cauchy-Schwartz, so that the integral in $\dTwo$ is inside the square root). The integral in the first term is bounded by the conformal charge as given in \eqref{eLabelForAngledConformalIn}. The second angular derivatives of $\MaxF$ will satisfy the same estimate, and we can use the second angular derivatives to control the value of the component, through a Sobolev estimate. Hence,
\begin{align*}
|\CoordWt(t,\rs,\theta,\phi)|
\leq& C\left(\int_{\{t\}\times\{\rs\}\times S^2} \sum_{k=0}^2 |\Lie_\SymGen^k\CoordWt|^2 \dTwo\right)^\frac12 \\
\leq& \left(\frac{\uout-\uin}{\uout\uin}\right)^{1/2} \left(\normKangN{4} +\normTangN{8}\right)^{1/2} ,\\
|\MyF(t,\rs,\theta,\phi)| \leq& \left(\frac{\uout-\uin}{\uout\uin} \right)^{1/2} r^{-2} \left(\normKangN{4} +\normTangN{8}\right)^{1/2}.
\end{align*}
Outside the light-cone, where $t<\rs$, we integrate over the curve $c_4$, and the end point value is replaced by $|\CoordWt(0,t+\rs,\theta,\phi)|$, which decays like $\rs'^{-1/2}=(t+\rs)^{-1/2}$ (since $\MyF$ decays like $\rs'^{-5/2}$). Thus, we have
\begin{align*}
|\MyF(t,\rs,\theta,\phi)| \leq& \left(\frac{\uout-|\uin|}{\uout|\uin|}\right)^{-1/2} r^{-2} \left(\normKangN{4} +\normTangN{8} + \SupInitData \right)^{1/2}.
\end{align*} 
In the region where $|\uin|<1$, instead of using $\GenEnergy{\vecConf}$, we could have used $\GenEnergy{\vecT}$, which does not have a vanishing factor of $\uin^2$. Thus, we may replace $(\uout-|\uin|)/\uout|\uin|$ by $(\uout-|\uin|)/\uout(1+|\uin|)$. 

Now, we prove decay for $\MyFp$ by again integrating along ingoing, radial, null geodesics again. From any given point, we integrate along $c_4$ to the endpoint where $t=0$, 
\begin{align*}
|\CoordWtp(t,\rs,\theta,\phi)|
\leq& \int_{c_4} |\vecncoord\CoordWtp | d\uin +|\CoordWtp(0,t+\rs,\theta,\phi)|\\
\leq& \left(\int_{c_4} r^{-2}d\uin \right)^\frac12 \left(\int_{c_4} \uout^2 |\vecncoord\CoordWtp|^2 r^2 d\uin\right)^\frac12 \uout^{-1}+|\CoordWtp(0,t+\rs,\theta,\phi)| \\
\leq& \left(\int_{c_4} \uout^2 |\vecncoord\CoordWtp|^2 r^2 d\uin\right)^\frac12 \uout^{-1}r^{-1/2}+|\CoordWtp(0,t+\rs,\theta,\phi)| .
\end{align*}
The endpoint will be bounded by $\rs'^{-3/2}=(t+\rs)^{-3/2}$. We now integrate in the angular variables, differentiate in the angular directions, and apply a spherical Sobolev estimate to get
\begin{align*}
|\CoordWtp(t,\rs,\theta,\phi)|
\leq& C \uout^{-1}r^{-1/2} \left(\normKangN{4}+\normTangN{8}\right)^{\frac12} \\
&+C\uout^{-3/2}\SupInitData ,\\
|\MyFp|\leq& C \uout^{-1}r^{-3/2} \left(\normKangN{4}+\normTangN{8}+\SupInitData \right)^{\frac12} .
\end{align*}

Finally, for $\MyFm$, we integrate along outgoing, radial, null rays on which $\uin$ is constant. Inside the light-cone $t>\rs$, we take the curve $c_2$ from $(t,\rs,\theta,\phi)$ to $(t-\rs,0,\theta,\phi)$. The estimate is 
\begin{align*}
|\CoordWtm(t,\rs,\theta,\phi)|
\leq& \int_{c_2} |\veclcoord\CoordWtm | d\uout +|\CoordWtm(t-\rs,0,\theta,\phi)|\\
\leq& (\int_{c_2} r^{-2}d\uout )^\frac12 (\int_{c_2} \uin^2 |\veclcoord\CoordWtm|^2 r^2 d\uout)^\frac12  \uin^{-1}+|\CoordWtm(t-\rs,0,\theta,\phi)|\\
\leq& C \uin^{-1}  (\int_{c_2} \uin^2 |\veclcoord\CoordWtm|^2 r^2 d\uout)^\frac12 +|\CoordWtm(t-\rs,0,\theta,\phi)| .
\end{align*}
In the stationary region, the decay rate is also $(t')^{-1}=(t-\rs)^{-1}$, so the decay rate is
\begin{align*}
|\CoordWtm(t,\rs,\theta,\phi)| \leq& C |\uin|^{-1} \left(\normKangN{4}+\normTangN{8}\right)^\frac12 , \\
|\MyFm|\leq& C |\uin|^{-1} r^{-1} \left(\normKangN{4}+\normTangN{8}\right)^\frac12 . 
\end{align*}
For $t<\rs$, a similar argument can be made by integrating along $c_3$, with the value at the other end point being $\CoordWtm(0,\rs-t,\theta,\phi)$, where we have faster decay, 
\begin{align*}
|\CoordWtm(t,\rs,\theta,\phi)| \leq& C |\uin|^{-3/2} \left(\normKangN{4}+\normTangN{8}\right)^\frac12 +|\uin|^{-3/2} \SupInitData, \\
|\MyFm|\leq& C |\uin|^{-3/2} r^{-1} \left(\normKangN{4}+\normTangN{8} +\SupInitData\right)^\frac12 . 
\end{align*}
Again, in the region $|\uin|<1$, we can use $\GenEnergy{\vecT}$ instead of $\GenEnergy{\vecConf}$ to get a better bound when $\uin$ vanishes. 
\end{proof}
\end{lemma}

We now turn to proving decay in the ``near'' region, $\rs<0$. Since for any fixed interval $2M<r_1<r<r_2$, we can apply the results from section \ref{sStationaryDecay}, the main purpose of the following lemma is to prove estimates which are uniform in $r$ so that they can be extended to the event horizon. Note that $\uout$ extends smoothly to the event horizon and is an affine parameter for tangential geodesics. 

Since the vector fields $\vecT$ and $\vecConf$ vanish on the bifurcation sphere, the boundedness of the associated energy allows rapid divergence of the (normalized) energy density near there. Not surprisingly, this is not sufficient to control the Maxwell field. In light of this, it is somewhat surprising that the energies associated with $\vecT$ and $\vecConf$ are sufficient to prove decay for the correctly normalized components of the Maxwell field tensor corresponding to $\MyFp$. 

As explained in the introduction, the correctly normalized basis for stating results on or near the event horizon is
\begin{align*}
\dt+\dr, &
&\horifac^{-1}(\dt-\dr), &
&r^{-1} \vecea, &
&r^{-1} \veceb .
\end{align*}
We can equally well replace $r^{-1}\vecea$ and $r^{-1}\veceb$ by $\vecm$ and $\vecmb$ or by $\vecHunit$ or $\vecPhunit$. For large $\uout$, this is the ``correctly normalized'' basis, because it is the result of parallelly transporting the original, normalized basis on the initial data surface, $t=0$, along ingoing null geodesics to reach the event horizon. 

The method used in the previous lemma gives decay rates of $\uout^{-1}$, $\uout^{-1/2}$, and $\uout^{-1}$ for $\CoordWtp$, $\CoordWt$, and $\CoordWtm$ respectively. The functions $\CoordWtp$ and $\CoordWt$ are correctly normalised (except for bounded factors) as $r\rightarrow2M$, but $\horifac^{-1}\CoordWtm$ is the correctly normalised component in this region. For $\CoordWt$, we prove a different preliminary decay rate and then use the divergence theorem to obtain a rate of $\uout^{-1}$. For the correctly normalised, negative-weight component, we use the $\uout^{-1}$ decay for $\CoordWt$, a transport equation, and an integrating factor to get $\uout^{-1}$ decay. We note that the vector field $\horifac^{-1}(\vecT-\vecR)$, which is a smoothed version of the vector field $Y$ in \cite{DafermosRodnianski} can be used to prove boundedness for this component without using a transport equation. 

\begin{lemma}[Decay for $\rs<0$.]
\label{lNearDecay}
There is a constant $C$ such that if $\MaxF$ is a solution of the Maxwell equations \eqref{eMaxwellEquationDiv}-\eqref{eMaxwellEquationAlt}, then for all $t\geq0$, $\rs<0$, $(\theta,\phi)\in S^2$ such that $\uout>1$, 
\begin{align*}
|\MaxF(\dt+\dr, \vecHunit)| +|\MaxF(\dt+\dr,\vecPhunit)| 
\leq& C\uout^{-1} \left(\normKangN{4}+\normTangN{8}\right)^\frac12, \\
|\MyF(t,\rs)|
\leq& C|\MaxF(\dt+\dr,\horifac^{-1}(\dt-\dr))| +|\MaxF(\vecHunit,\vecPhunit)| \\
\leq& C\uout^{-1}\left(\normKangN{4}+\normTangN{8} +\normTunitBifSphN{3}\right)^\frac12 , \\
F(\horifac^{-1}(\dt-\dr),\vecHunit +i\vecPhunit) \leq& C \uout^{-1} \left(\normKangN{4} + \normTangN{8} + \normTunitBifSphN{3} \right)^{1/2} .
\end{align*}
\begin{proof}
The first part of this proof is similar to that of lemma \ref{lFarDecay}. The main difference is that we must track factors of $\horifac$ carefully, but we may ignore factors of $r$ since it is bounded above and below by positive constants. When tracking factors of $\horifac$, we use $\horifac\MyFi(t',\rs',\theta',\phi')$ to denote the value of $\horifac\MyFi$ at $(t',\rs'\theta',\phi')$ even if an unprimed set of coordinates is in use simultaneously. Since we are only considering $\uout>1$, the ingoing, radial, null rays from any point will hit the stationary region $\rs=0$. This is illustrated in figure \ref{FigNearDecaypm}

\begin{figure}
\begin{center}
\input{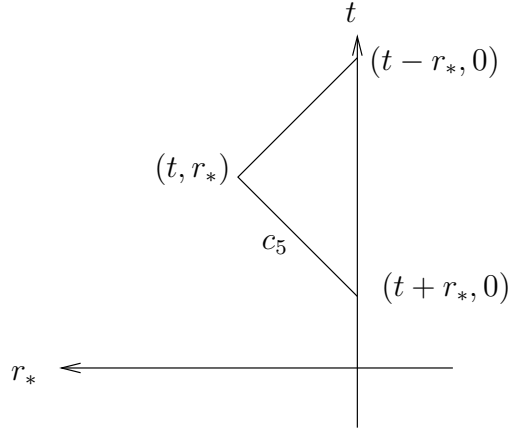}
\caption{Null rays in the inner region $\rs<0$. The angular variables have been suppressed. The curve $c_5$ goes from a point in the stationary region $\rs=0$ to an arbitrary point in the regions $t>0$, $\rs<0$, $\uout>1$ along an ingoing, null, radial geodesics. }
\label{FigNearDecaypm}
\end{center}
\end{figure}

Integrating along an ingoing, radial, null geodesic, $c_5$, from $(t,\rs,\theta,\phi)$ to $(t+\rs,0,\theta,\phi)$, in the same way as in the proof of lemma \ref{lFarDecay}, we have
\begin{align*}
|\CoordWtp(t,\rs,\theta,\phi)|
\leq& \int_{c_5} |\vecncoord\CoordWtp| d\uin +|\CoordWtp(t+\rs,0,\theta,\phi)|\\
\leq& \left(\int_{c_5} \horifac d\uin\right)^\frac12 
\left(\int_{c_5} |\vecncoord\CoordWtp|^2 \uout^2 r^2 \horifac^{-1} d\uin\right)^\frac12 \uout^{-1} +|\CoordWtp(t+\rs,0,\theta,\phi)|  .
\end{align*}
The integral of $\horifac$ with respect to $d\uin$ is, up to a factor of $2$, the same as the integral of $\horifac$ with respect to $d\rs$. By a change of variables, this is the integral of $1$ with respect to $dr$. Thus, the contribution from the first integral in the first term is bounded by a constant. After integrating in the angular variables, the second integral is bounded by the conformal energy of the angular derivatives of $\MaxF$. The second term is the value of the component in the stationary region, so it decays like $(t+\rs)^{-1}$. Applying the angular derivative and Sobolev estimate argument from the previous lemma, 
\begin{align*}
|\CoordWtp|
=|\MaxF(\dt+\dr, \vecHunit+i \vecPhunit)|
\leq& C \uout^{-1} \left(\normKangN{4}+\normTangN{8}\right)^\frac12 .
\end{align*}
This proves the first result. 

For the zero weight component, we first prove an intermediate result for $\horifac^{1/2}\MyF$. From integrating along a surface of constant $t$, we have
\begin{align*}
\horifac^{1/2}|\uout|^2 \MyF^2(t,\rs,\theta,\phi)
=& -\int_{\{t\}\times[\rs,0]\times\{\theta\}\times\{\phi\}} \dr\left(\horifac^{1/2}|\uout|^2 \MyF^2\right) d\rs'\\
& + \horifac^{1/2}\uout^2 \MyF^2(t+\rs,0,\theta,\phi). 
\end{align*}
The integrand can be estimated by dropping negative terms and applying the Cauchy-Schwartz inequality as
\begin{align*}
-\dr\left(\horifac^{1/2}\uout^2|\MyF|^2\right)
=& -\frac12 \horifac^{1/2}\frac{2M}{r^2} \uout^2 |\MyF|^2 
- 2\horifac^{1/2}\uout |\MyF|^2 \\
&- 2\horifac^{1/2} \uout^2 \Re(\MyF \dr\MyF) \\
\leq& C\left( \uout^2 \horifac |\MyF|^2 + \uout^2 |\dr\MyF|^2 \right) .
\end{align*}
Thus, for $\uout>1$, the integrand is controlled by the conformal energy using the standard Sobolev estimate argument by
\begin{align*}
-\int_{\{t\}\times[\rs,0]\times\{\theta\}\times\{\phi\}} \dr\left(\horifac^{1/2}|\uout|^2 \MyF^2\right) d\rs'
\leq&C\left( \normKangN{4} + \normTangN{8} \right).
\end{align*}
The end point can be controlled by the stationary decay result, 
\begin{align*}
\horifac^{1/2}\uout^2|\CoordWt(t+\rs,0,\theta,\phi)|^2 \leq& C \left(\normKangN{4} +\normTangN{8} \right), 
\end{align*}
so that
\begin{align*}
\MyF(t,\rs,\theta,\phi)
\leq& C\horifac^{-1/4} \uout^{-1} \left( \normKangN{4} + \normTangN{8} \right)^{1/2}.
\end{align*}

\begin{figure}
\begin{center}
\input{FigDivThmRegion.pstex_t}
\caption{The region $\DivThmRegion$.}
\end{center}
\label{FigDivThmRegion}
\end{figure}

This estimate can now be improved. For a given point $(t,\rs,\theta,\phi)$, consider the two-dimensional surface
\begin{align*}
\DivThmRegion=&\{(t',\rs',\theta,\phi) : t'\geq0,t'+\rs'\leq t+\rs, t-\rs\leq t'-\rs'\leq t+\rs \} . 
\end{align*}
This is illustrated in figure \ref{FigDivThmRegion}. Applying the (two-dimensional) divergence theorem with the vector field $(\vecncoord\CoordWt)\vecncoord$, we have
\begin{align*}
\CoordWt(t,\rs,\theta,\phi) -\CoordWt(t+\rs,0,\theta,\phi) +\int_{\{0\}\times[\rs-t,-t-\rs]\times\{\theta\}\times\{\phi\}} \vecncoord\CoordWt d\rs 
=& -2 \int_\DivThmRegion \veclcoord\vecncoord\CoordWt d\rs dt .  
\end{align*}
From the wave equation \eqref{ePriceWaveEqn} for $\CoordWt$, we have
\begin{align*}
\CoordWt(t,\rs,\theta,\phi)
=& \CoordWt(t+\rs,0,\theta,\phi)
+\int_{\{0\}\times[\rs-t,-t-\rs]\times\{\theta\}\times\{\phi\}} \vecncoord\CoordWt d\rs\\
&-2 \int_\DivThmRegion r^{-2}\horifac(-\slap)\CoordWt d\rs dt .
\end{align*}

We now estimate the terms on the right. The first is bounded by the stationary decay result. The second is an integral in the initial data surface $t=0$, so it can be controlled by integrals of the initial data. We have
\begin{align*}
|\int_{[\rs-t,-t-\rs]} \vecncoord\MyF d\rs|
\leq& \left(\int_{[-\infty,-t-\rs]} \horifac^{1/2} d\rs\right)^{1/2} \left(\int_{[-\infty,-t-\rs]} \horifac^{-1/2}|\vecncoord\MyF|^2 d\rs\right)^{1/2} . 
\end{align*}
The first integral is bounded by 
\begin{align*}
\int_{[-\infty,-t-\rs]} \horifac^{1/2} d\rs
\leq& C\horifac^{1/2} \leq C e^{(-t-\rs)/2M} \leq C \uout^{-2} .
\end{align*}
Since $\horifac^{-1/2}\vecncoord$ and $\partial/\partial\Buin$ differ only by smooth functions of $r$, the second can be estimated by 
\begin{align*}
\int_{[-\infty,-t-\rs]} \horifac^{-1/2}|\vecncoord\MyF|^2 d\rs
\leq& C \int_{[-\infty,-t-\rs]} |\frac{\partial}{\partial\Buin} \MyF|^2 \horifac^{1/2} d\rs . 
\end{align*}
Since the stationary tetrad and the one based on the $\Buout, \Buin, \theta, \phi$ coordinate system differ only by smooth functions of $r$, 
\begin{align*}
\int_{[-\infty,-t-\rs]} \horifac^{-1/2}|\vecncoord\MyF|^2 d\rs
\leq& C \int_{[-\infty,-t-\rs]} |\frac{\partial}{\partial\Buin}\left(\frac12\MaxF(\veclunit,\vecnunit)+\MaxF(\vecHunit,\vecPhunit)\right) |^2  \horifac^{1/2}d\rs \\
\leq& C \int_{[-\infty,-t-\rs]} |\frac{\partial}{\partial\Buin}\left(\frac12\MaxF(\frac{\partial}{\partial\Buout},\frac{\partial}{\partial\Buin})+\MaxF(\vecH,\frac{1}{\sin\theta}\vecPh)\right) |^2 \horifac^{1/2} d\rs \\
&+C \int_{[-\infty,-t-\rs]} |\frac12\MaxF(\frac{\partial}{\partial\Buout},\frac{\partial}{\partial\Buin})+\MaxF(\vecH,\frac{1}{\sin\theta}\vecPh) |^2 \horifac^{1/2} d\rs .
\end{align*}
Since coordinate vector fields commute, 
\begin{align*}
\int_{[-\infty,-t-\rs]} \horifac^{-1/2}|\vecncoord\MyF|^2 d\rs
\leq& C \int_{[-\infty,-t-\rs]} \sum_{k=0}^1 |\Lie_{\frac{\partial}{\partial\Buin}}^{k} \MaxF|_{\BifurcationCoordinateBasis}^2 \horifac^{1/2} d\rs .
\end{align*}
The same argument could have been applied to the second angular derivatives of $\MaxF$, which could have been used in a Sobolev estimate. This would have lead to
\begin{align*}
\int_{[-\infty,-t-\rs]} \horifac^{-1/2}|\vecncoord\MyF|^2 d\rs
\leq& C \normTunitBifSphN{3} .
\end{align*}
Thus, the integral along the initial time slice is bounded by
\begin{align*}
|\int_{[\rs-t,-t-\rs]} \vecncoord\MyF d\rs|
\leq& C \uout^{-1} \left( \normTunitBifSphN{3} \right)^{1/2} .
\end{align*}
Finally, we estimate the integral over $\DivThmRegion$ by breaking it into two parts, $\DivThmRegionA=\DivThmRegion\cap \{t>2|\rs|\}$ and $\DivThmRegionB=\DivThmRegion\cap \{t\leq2|\rs|\}$. In $\DivThmRegion$, 
\begin{align*}
|\int_\DivThmRegion& r^{-2}\horifac(-\slap)\CoordWt d\rs dt| \\
\leq& \sup_{\DivThmRegionA}\left(\horifac^{1/4} (-\slap)\CoordWt\right) \int_\DivThmRegionA r^{-2}\horifac^{3/4} d\rs dt \\
&+ \left( \int_\DivThmRegionB r^{-2}\horifac d\rs dt\right) \left(\int_\DivThmRegionB r^{-2}\horifac|\slap\CoordWt|^2 d\rs dt  \right)
\end{align*}
On the first line of the right-hand side, the supremum term decays like $\uout^{-1}$ by the intermediate result, and the integral term is uniformly bounded by the exponential decay of $\horifac$ with respect to $\rs$. In the second line, the second integral is bounded by estimate \eqref{eWaveLocalDecay}, and the first integral is bounded by $\horifac$ evaluated at the point $(t',\rs')=((2/3)\uout,(1/3)\uout)$, and hence decays faster than $\uout^{-1}$. Combining all these results gives
\begin{align*}
\MyF
\leq& C \uout^{-1} \left(\normKangN{4} +\normTangN{8} + \normTunitBifSphN{3} \right) .
\end{align*}

We begin our analysis of $\CoordWtm$ with an intermediate decay result, using the same sort of simple argument as was used for $\CoordWtp$. Integrating along an outgoing, radial, null ray, we have
\begin{align*}
\CoordWtm(t,\rs,\theta,\phi)
=& \int |\veclcoord\CoordWtm| d\uout + \CoordWtm(t-\rs,0,\theta,\phi) \\
|\CoordWtm(t,\rs,\theta,\phi)|
\leq& \uin^{-1} \left(\int \horifac d\uout\right)^{1/2} \left( \int |\veclcoord\CoordWtm|^2 \uin^2 \horifac^{-1} d\uout\right)^{-1} + |\CoordWtm(t-\rs,0,\theta,\phi)| \\
\leq& C \uin^{-1} \left(\normKangN{4} + \normTangN{8}\right)^{1/2} . 
\end{align*}
Since we are working in the inner region, $\rs\leq0$, there is the estimate $\uout<\uin$, and 
\begin{align*}
\CoordWtm(t,\rs,\theta,\phi)
\leq& C \uout^{-1} \left(\normKangN{4} + \normTangN{8}\right)^{1/2} . 
\end{align*}

To obtain stronger estimates, we will need to integrate along outgoing null geodesics starting near the initial data surface. We work with the $\Buout$ and $\Buin$ coordinates to control the correctly normalised, negative-weight component near the bifurcation sphere. If bounded initial data is posed on the surface $\Buout\Buin=-1$, which corresponds to the union of the $t=0$ surfaces in the two exterior regions, then, at least in in some small neighborhood of the bifurcation sphere, the components of the Maxwell field with respect to the smooth $(\Buout,\Buin,\theta,\phi)$ coordinate system must remain bounded by a multiple of their initial value. This is essentially a Cauchy stability result, as was used for the wave equation \cite{DafermosRodnianski}. Thus, in some sufficiently small neighborhood of the bifurcation sphere, in the exterior, the outgoing component is bounded by
\begin{align*}
4M e^{u-/4M} |\CoordWtm|
\leq&C_1|\MaxF(\frac{\partial}{\partial \Buin}, \vecH)|
+|\MaxF(\frac{\partial}{\partial \Buin}, \vecPh)| \nonumber\\
\leq& C\normTunitBifSphN{3} .
\end{align*}
In particular, we can pick a ${\uout}_0\ll0$ such that on the hypersurface $\uout={\uout}_0$ where $t\geq0$, 
\begin{align}
\horifac^{-1} |\CoordWtm|
\leq& C\normTunitBifSphN{3} .
\label{ePrelimBndOnMyFmNearHorizon}
\end{align}

We now use the decay for $\CoordWt$ and $\CoordWtm$ to prove a stronger estimate. From the Price equation \eqref{ePriceiv}, we have
\begin{align*}
\veclcoord\left( e^{\uout/4M}\horifac^{-1}\CoordWtm\right)
=& \left(\frac{1}{2M}-\frac{2M}{r^2}\right)\left( e^{\uout/4M}\horifac^{-1}\CoordWtm\right)
+ e^{\uout/4M} r^{-2} \CoordWt . 
\end{align*}
Since $1/2M -2M/r^2$ vanishes linearly at $r=2M$, it is bounded by $C\horifac$. Integrating along an outgoing geodesic starting on $\uout={\uout}_0$ at $(t_0,{\rs}_0,\theta,\phi)$ and going to $(t,\rs,\theta,\phi)$, we have
\begin{align*}
|e^{\uout/4M}\horifac^{-1}\CoordWtm(t,\rs,\theta,\phi)|
\leq& C\int_{c_6} e^{\uout/4M} |\CoordWtm+\CoordWt| d\uout + |e^{{\uout}_0/4M}\horifac^{-1}\CoordWtm(t_0,\rs_0,\theta,\phi)| . 
\end{align*}
The geodesic along which we integrate is illustrated in figure \ref{FigWtm}. We break the integral into two pieces, with one going from ${\uout}_0$ to $\uout/2$ and the other going from $\uout/2$ to $\uout$. From the boundedness and decay for $\CoordWt$ and $\CoordWtm$, the integral is bounded by
\begin{align*}
\int_{c_6} e^{\uout/4M} |\CoordWtm+\CoordWt| d\uout
\leq& C (e^{\uout/8M} + \uout^{-1} e^{\uout/4M}) \left(\normKangN{4}+\normTangN{8} +\normTunitBifSphN{3}\right)^{1/2} .
\end{align*}
By estimate \eqref{ePrelimBndOnMyFmNearHorizon}, the end point, $\horifac^{-1}\CoordWtm(t_0,{\rs}_0,\theta,\phi)$ is bounded. Thus, 
\begin{align*}
|\horifac^{-1}\CoordWtm(t,\rs,\theta,\phi)|
\leq&C (e^{-\uout/8M} + \uout^{-1} ) \left(\normKangN{4}+\normTangN{8} +\normTunitBifSphN{3}\right)^{1/2} \\
&+ C e^{(-\uout+{\uout}_0)/4M} \left(\normKangN{4}+\normTangN{8} +\normTunitBifSphN{3}\right)^{1/2} \\
\leq& C \uout^{-1}\left(\normKangN{4}+\normTangN{8} +\normTunitBifSphN{3}\right)^{1/2} .
\end{align*}
This provides the desired decay. 

\begin{figure}
\begin{center}
\input{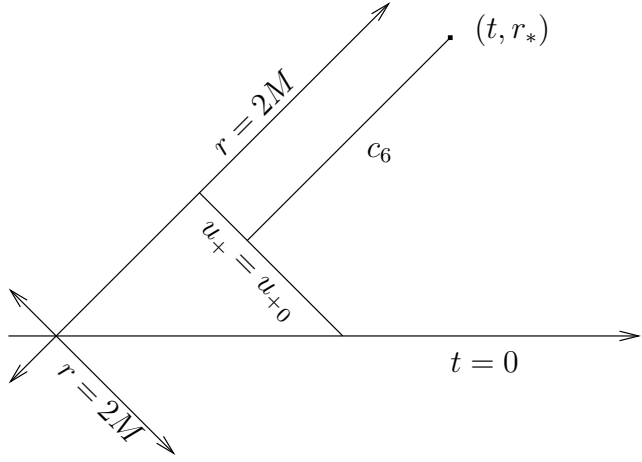}
\caption{An outgoing null ray from the hypersurface $\uout={\uout}_0$ to the point under consideration. The angular variables have been suppressed. The region under consideration contains a portion near the bifurcation sphere and is near the event horizon. Therefore, the $t$ and $\rs$ coordinates are not used. Instead the hypersurfaces $t=0$, $r=2M$, and $\uout={\uout}_0$ are indicated. }
\label{FigWtm}
\end{center}
\end{figure}
\end{proof}
\end{lemma}

Combining the two lemmas in this section, we have theorem \ref{tNearAndFarDecay}. 

\begin{remark}
\label{SimplifiedInitialData}
Finally, we provide a simpler sufficient condition for the initial data to have finite norm. 

The energies $\GenEnergy{\vecT}[\MaxF]$ and $\GenEnergy{\vecConf}[\MaxF]$ are the integrals of the field components squared, $\MyFi^2$, against the weight $\horifac r^2$ and $\rs^2 \horifac r^2$ respectively. Thus, if the field components are bounded everywhere and decay at least as fast as $r^{-(5/2+\epsilon)}$ (for $\epsilon>0$) as $r\rightarrow\infty$, these energies will be finite. Thus, if the Maxwell field and its first eight derivatives decay at this rate, then 
\begin{align*}
\normKangN{4} + \normTangN{8} 
\end{align*}
will be bounded, and the result of theorem \ref{tDecayInStationaryRegionsInText} will apply. 

Similarly, if the field components decay at least as fast as $r^{-(5/2+\epsilon)}$, then $\SupInitData$ will be trivially bounded, so that the results of lemma \ref{lFarDecay} will give decay in the far region $\rs>1$. 

The energy $\GenEnergy{\vecTunit}$ is the integral of the field components squared, $\MyFi^2$, against the measure $\horifac^{1/2}r^2 d\rs$. Thus, if the components are bounded as $\rs\rightarrow-\infty$, then this energy will be finite. Since (for $\rs<0$), the vector fields in $\BifurcationCoordinateBasis$ are coordinate vector fields extending in a neighbourhood of the bifurcation sphere, they are smooth vector fields. Thus, if $\MaxF$ and its first three derivatives with respect to this collection of smooth vector fields have finite components, then 
\begin{align*}
\normTunitBifSphN{3}
\end{align*}
will be finite, and the results of lemma \ref{lNearDecay} will apply. 

This verifies the footnote to theorem \ref{tDecayInStationaryRegions} that boundedness and $r^{-5/2+\epsilon}$ decay for $\MaxF$ and its first eight derivatives on the initial data surface $t=0$ is sufficient to prove the decay results in this paper. The same remark applies to the result in theorem \ref{tNearAndFarDecay}. 
\end{remark}

%
%
\appendix
\section{Exclusion of the non-radiatable mode of the Maxwell field}
\label{sExclusionOfNonRadiatable}
In this section, we show that if the Maxwell field has finite conformal charge, it has no spherically symmetric part and explain why this is a physically reasonable assumption. The absence of dynamic, spherically symmetric components is well-known in the literature \cite{Price}. In $\Reals^{1+3}$, since the electric and magnetic fields are divergence free, the spherically symmetric component of the Maxwell field is always zero. While there are spherically symmetric solutions on the Schwarzschild manifold, we show that they have no dynamics, since these solutions are constant in $t$. There is a two-parameter family of such solutions described by the central electric and magnetic charge. These solutions do not vanish on the event horizon and decay like $1/r^2$ at infinity, so they fail to be in the finite conformal energy class we consider. 

The Maxwell field $\MaxF$ can be written as
\begin{align*}
\MaxF(t,\rs,\theta,\phi)
=& \Myrho(t,\rs,\theta,\phi) \horifac \Extd t\wedge \Extd \rs\\
&+ r\horifac^{1/2}\Extd t\wedge\omega_0(t,\rs,\theta,\phi)\\
&+ r\horifac^{1/2}\Extd \rs\wedge\omega_1(t,\rs,\theta,\phi)\\
&+\Mysigma(t,\rs,\theta,\phi) r^2 \Omega ,\\
*\MaxF(t,\rs,\theta,\phi)
=& -\Mysigma \horifac \Extd t\wedge \Extd \rs\\
&- r\horifac^{1/2}\Extd t\wedge\SHD\omega_1\\
&+ r\horifac^{1/2}\Extd \rs\wedge\SHD\omega_0\\
&+ \Myrho r^2 \Omega ,
\end{align*}
with $\Myrho$ and $\Mysigma$ scalar functions, with $\SHD$ the Hodge dual on $S^2$, with $\omega_0$ and $\omega_1$ in $\forms^1(S^2)$ for each value of $t$ and $\rs$, and with $\Omega$ the standard volume form on $S^2$. We have used a stationary, instead of null, decomposition, so $\omega_0$ and $\omega_1$ appear instead of $\CKalpha$ and $\CKalphab$. The scalars $\CKrho$ and $\CKsigma$ are the standard ones from the null decomposition. 

We first remind the reader that there is no $\omega\in\forms^1(S^2)$ with $\ExtdS\omega= \SHD C_0$. If there were one, we could write this condition in coordinates: 
\begin{align*}
\omega=&\omega_\theta \Extd\theta +\omega_\phi\Extd\phi ,\\
\ExtdS \omega=&\SHD C_0 ,\\
\omega_{\theta,\phi} - \omega_{\phi,\theta} 
=& C_0 \sin(\theta) .
\end{align*}
Let $f(\theta)=\int_{\phi\in S^1} \omega_\phi \Extd \phi$. Since $\omega$ is smooth, $f$ is continuous on $[0,2\pi]$ and vanishing at $0$ and $\pi$ (since the integral is over a single point in $S^2$ in these cases). From $\Extd \omega=\SHD C_0$, we have $f'(\theta)=C_0 \sin(\theta)$, and $f(\pi)-f(0)=\int_{0}^{\pi} C_0 \sin(\theta) \Extd \theta >0$. Thus, the condition that $f$ vanishes at $0$ and $\pi$ implies $C_0=0$. 

The Maxwell equations are (taking the components orthogonal to various $1$-forms)
\begin{align*}
\text{Equation }&\text{ Orthogonal $1$-form}\\
\Extd F=0:
&\Extd t:& 0=&r\horifac^{1/2}\Extd\rs\wedge \ExtdS\omega_1 + (\dr(\Mysigma r^2)) \Extd\rs\wedge\Omega  \\
&\Extd \rs:& 0=&r\horifac^{1/2}\Extd t\wedge \ExtdS\omega_0 + (\dt(\Mysigma r^2)) \Extd t\wedge\Omega  \\
\Extd *F=0:
&\Extd t:& 0=&r\horifac^{1/2}\Extd\rs\wedge \ExtdS\SHD\omega_0 + (\dr(\Myrho r^2)) \Extd\rs\wedge\Omega  \\
&\Extd \rs:& 0=&r\horifac^{1/2}\Extd t\wedge \ExtdS\SHD\omega_1 + (\dt(\Myrho r^2)) \Extd t\wedge\Omega
\end{align*}
Since $\Extd t$ and $\Extd\rs$ are spherically symmetric, these can be projected onto the $l=0$ spherical harmonic (equivalently, we can contract on the $S^2$ volume). Since there is no $l=0$ component for the $1$-forms $\omega_0$, $\omega_1$, $\SHD\omega_0$ and $\SHD\omega_1$, we find, 
\begin{align*}
\dr(\Mysigma_{l=0} r^2)=& 0\\
\dt(\Mysigma_{l=0} r^2)=&0 \\
\dr(\Myrho_{l=0} r^2)=&0 \\
\dt(\Myrho_{l=0} r^2)=&0 .
\end{align*}
From which it follows that the $l=0$ components are given by integration constants $q_E$ and $q_B$, 
\begin{align*}
\Myrho_{l=0}=& \frac{q_{E}}{r^2} \\
\Mysigma_{l=0}=& \frac{q_{B}}{r^2} .
\end{align*}
Thus, there is no dynamics is the $\SHarmPara=0$ mode, since the $t$ derivative is always zero. These solutions do not decay sufficiently rapidly to have finite conformal energy. 

The exclusion of these spherically symmetric solutions is physically reasonable. Physically, the solutions represent a perturbation of the Schwarzschild black hole to a charged Reissner-Nordstrom solution, not an external perturbation by radiation. Price refers to these spherically symmetric solutions as the ``non-radiatable'' modes, since the solutions in this two parameter family are static. Since the Maxwell equations are linear and commute with angular derivatives, the spherically symmetric component does not couple to the other components, so it will not affect the rest of our analysis to eliminate the spherically symmetric components. In analogy with the theory of solitons, we might think of the Reissner-Nordstrom solutions as a manifold in the space of solutions to the Maxwell-Einstein system. In this case, the decoupled Maxwell equations with $\SHarmPara>0$ correspond to linearized perturbations from this manifold, whereas perturbations with $\SHarmPara=0$ correspond to linearized perturbations along the manifold of stationary solutions. 

\section{Analysis of the wave equation}
\label{sOneDWaveAnalysis}
We now prove decay estimates for solutions to the wave equation \eqref{ePriceWaveEqn}. For this equation, there is also an energy and conformal energy, which we use in our analysis. As with the Maxwell field, we must control the trapping of $\SFn$ near the photon sphere to control the growth of the conformal charge. We do this with a local decay estimate and employ light-cone localization to obtain a local decay estimate of the full strength we require. The arguments and results of this section are a slight modification of the argument in \cite{BlueSterbenz}, only, in this case, the situation is simpler. 

For the wave equation on the Schwarzschild manifold, previous analysis \cite{BlueSofferLongPaper,BlueSterbenz,DafermosRodnianski} has required a decomposition onto spherical harmonics. On each spherical harmonic, the wave equation can be treated as a one-dimensional wave equation with an effective potential. The main estimate uses a vector field, $\LDmult$, which points away from the maximum of the effective potential. In the case of the geometrically defined wave equation, the location of these maxima depend on the spherical harmonic parameter, and $\LDmult$ has been modified to fit each spherical harmonic. The equation \eqref{ePriceWaveEqn} is simpler, with the maxima of the effective potential always at $\rs=0$. Thus, a very minor modification of the previous analysis allows us to make the estimate without using a spherical harmonic decomposition. 

Since the potential $V_L= \frac{1}{r^2}\horifac$ is real-valued, we may analyse the real and complex parts, which each satisfy \eqref{ePriceWaveEqn}, separately. Thus, we may assume our solutions are real-valued. 

We use the method of multipliers to analyse \eqref{ePriceWaveEqn}. Although it maybe possible to introduce a Lagrangian formulation and an energy-momentum tensor, we do not do so because this is not the geometrically defined wave equation $\CoDeriv^\ia\CoDeriv_\ia(r^{-1}u)$ so the Lagrangian for this system is quite artificial, because there would be a confusion between the energy-momentum for $\SFn$ and that of the full Maxwell field, and because some of the energies would require correction terms. The essence of the method is to choose a ``multiplier'' (a differential operator), apply it to the function $\SFn$, multiply by the equation, and integrate by parts. The most useful differential operators are typically those given by the vector fields from the Lagrangian method. 

We begin by recalling the energy and conformal energy, which were defined in section \ref{sEnergies}. Conservation of energy follows simply from the method of multipliers with the multiplier $\vecT=\dt$. Multiplying \eqref{ePriceWaveEqn} by $\dt\SFn$, we find
\begin{align}
0
=& (\dt\SFn)(-\dtt \SFn +\drr\SFn +\frac{1}{r^2}\horifac\Delta_{S^2}\SFn) \nonumber \\
=& -\frac12 \dt( |\dt\SFn|^2 +|\dr\SFn|^2 +V_L|\dAng\SFn|^2 )
+\dr(\dt\SFn\dr\SFn) 
+\dAng\cdot(V_L\dt\SFn\dAng\SFn) \label{eEnergyDivForm}.
\end{align}
 Integrating over a space-time slab gives conservation of energy: 
\begin{align*}
\PriceEnergy[\SFn](t)-\PriceEnergy[\SFn](0)=&0 .
\end{align*}
Using the method of multipliers with $\vecConf$ gives
\begin{align}
\PriceConf[\SFn](t_2)-\PriceConf[\SFn](t_1)
=&\intTslab{t_1}{t_2} 2t(2V_L+\rs V_L')|\dAng\SFn|^2 \dFourPOneD \nonumber\\
\leq&\intTslab{t_1}{t_2} t\chiTrap|\dAng\SFn|^2 \dFourPOneD  .\label{ePriceConformalIdentity}
\end{align}
This is similar to the estimate for the Maxwell equations. For this analysis, it is useful to introduce an energy localized inside the light cone. We let
\begin{align*}
\PriceMinEnergy=& \intTLCslice{t} \PriceEnergyDensity  \dThreePOneD .
\end{align*}

We need a variety of Hardy estimates. 

\begin{lemma}
If  $t\geq1$, $\chiHardy$ is a non-negative function which is positive in some open set $|\rs|<t$, and $\alpha>0$, then if $\HardyFn:\Reals\times S^2$ is a smooth function, and $\SFn:[t_1,t_2]\times\Reals\times S^2\rightarrow\Reals$ is smooth with $t\in[t_1,t_2]$ and $\SFn(t)=\HardyFn$, then 
\begin{align*}
\intTLCslice{t} \frac{|\HardyFn|^2}{(1+\rs^{2})} \dThreePOneD
\leq C\PriceMinEnergy[\SFn](t) , \\
\intTLCtighterslice{t} \frac{|\HardyFn|^2}{(1+|\rs|)^{\alpha+2}} \dThreePOneD
\leq C \intTLCtighterslice{t} \frac{|\dr \HardyFn|^{2}}{(1+\rs)^{\alpha}} + \chiHardy|\HardyFn|^2 \dThreePOneD
\end{align*}
\begin{proof}
We start working with $\alpha\geq0$. When ${\rs}_1>0$, 
\begin{align*}
\frac{|\HardyFn({\rs}_1)|^2}{(1+{\rs}_1)^{\alpha+1}} - |\HardyFn(0)|^2 
=& \int_0^{{\rs}_1} \dr \frac{|\HardyFn|^2}{(1+\rs)^{\alpha+1}} d\rs \\
=&\int_0^{{\rs}_1} \frac{2\HardyFn\dr \HardyFn}{(1+\rs)^{\alpha+1} } -(\alpha+1)\frac{|\HardyFn|^2}{(1+\rs)^{\alpha+2} } d\rs\\
\leq& \frac{\alpha+1}{2}\int_0^{{\rs}_1} \frac{|\HardyFn|^2}{(1+\rs)^{\alpha+2}} d\rs 
+ \frac{2}{\alpha+1} \int_0^{{\rs}_1}\frac{|\dr \HardyFn|^2}{(1+\rs)^{\alpha}} d\rs \\
&-(\alpha+1) \int_0^{{\rs}_1} \frac{|\HardyFn|^2}{(1+\rs)^{\alpha+2}} d\rs\\
\int_0^{{\rs}_1} \frac{|\HardyFn|^2}{(1+\rs)^{\alpha+2}} d\rs
\leq& \frac{4}{(\alpha+1)^2} \int_0^{{\rs}_1}\frac{|\dr \HardyFn|^2}{(1+\rs)^{\alpha}} d\rs 
+\frac{2}{\alpha+1} |\HardyFn(0)|^2 .
\end{align*}
Since, for any exponent $\beta\geq0$, $(1+\rs)^{\beta}$ is equivalent to $(1+\rs^2)^{\beta/2}$ on $[0,\infty)$, the powers of $(1+\rs)$ can be replaced by $(1+\rs^2)^{1/2}$. By symmetry, the same result holds on $(-{\rs}_1,0]$. Since $(1+\rs^2)^{-\beta}$ is uniformly equivalent to $(1+(\rs-{\rs}_0)^2)^{-\beta}$ for ${\rs}_0$ in a finite interval, the $|\HardyFn(0)|^2$ term can be replaced by $|\HardyFn({\rs}_0)|^2$ in any fixed interval. 

For $\alpha=0$, we take ${\rs}_1=(3/4)t$. By integrating the estimate over ${\rs}_0$ with ${\rs}_0$ in $(1/2,3/4)$, where $V_L$ is strictly positive, and then integrating over the angular variables, we find 
\begin{align*}
\intTLCslice{t} \frac{|\HardyFn|^2}{(1+\rs^{2})} \dThreePOneD
\leq C\PriceMinEnergy[\SFn](t) .
\end{align*}
Similarly, taking ${\rs}_1=(1/2)t$ and $\alpha>0$, for any non-negative function, $\chiHardy$, which is positive in some open set inside $|\rs|\leq t$,  
\begin{align*}
\intTLCtighterslice{t} \frac{|\HardyFn|^2}{(1+\rs)^{\alpha+2}} \dThreePOneD
\leq C \intTLCtighterslice{t} \frac{|\dr \HardyFn|^2}{(1+|\rs|)^{\alpha}} + \chiHardy|\HardyFn|^2 \dThreePOneD
\end{align*}
\end{proof}
\end{lemma}

We now prove a local decay estimate to control the trapping terms. To do this, we use a radial multiplier $\LDmult$ in terms of a weight $\LDwt$, 
\begin{align*}
\LDmult= \LDwt \dr + (\dr \LDwt)/2 .
\end{align*}
Assuming that the weight $\LDwt$ is a function of the $t$ and $\rs$ variables only, we have, 
\begin{align*}
-\dt(2\dot{\SFn}\LDmult \SFn)
=& -\dt(2\dot{\SFn}\LDwt \dr \SFn +\dot{\SFn}(\dr\LDwt) \SFn) \\
=&-\dr\left(\LDwt (\dr \SFn)^2 +(\dr \SFn)(\dr \LDwt)\SFn -V_L \LDwt (\dAng \SFn)^2 +\LDwt(\dot{\SFn})^2 \right)
-\dAng\cdot(V_L(\dAng \SFn)\LDwt \SFn) \\
&+2(\dr\LDwt)(\dr \SFn)^2 -\frac{(\drrr\LDwt) \SFn^2}{2} -(\dr V_L)\LDwt|\dAng \SFn|^2\\
&-2\dot{\SFn}\dot{g}(\dr \SFn)-\dot{\SFn}(\dr\dot{g})\SFn .
\end{align*}
We use the notation 
\begin{align*}
\PriceLDEnergy[\SFn](t)=\intTslice{t} \dot{\SFn}(\LDmult\SFn) \dThreePOneD . 
\end{align*}
Integrating over a space-time slab, 
\begin{align}
-2\PriceLDEnergy[\SFn](t_2)+2\PriceLDEnergy[\SFn](t_1)
= \intTslab{t_1}{t_2} \big(2(\dr\LDwt)(\dr \SFn)^2 -\frac{(\drrr\LDwt) \SFn^2}{2} -(\dr V_L)\LDwt|\dAng \SFn|^2&\nonumber\\
-2\dot{\SFn}\dot{g}(\dr \SFn)-\dot{\SFn}(\dr\dot{g})\SFn\big) &\dFourPOneD .\label{eGeneralLDFormula}
\end{align}

Taking $\chiLCproto$ to be smooth, non-negative, compactly supported in $[-3/4,3/4]$, and identically $1$ on $[-1/2,1/2]$, $\LDpara$ to be a sufficiently small parameter to be chosen later, and $\LDexp\in(1,2]$, we set
\begin{align*}
\chiLC=& \chiLCproto(\frac{\rs}{t}) ,\\
\LDsubwt=& \int_0^{\rs} \frac{1}{1+b|y|^\LDexp} dy ,\\
\LDwt(t,\rs)=& t \LDsubwt\chiLC .
\end{align*}
We now expand and group the terms on the right of \eqref{eGeneralLDFormula}, 
\begin{align}
-2\PriceLDEnergy[\SFn](t_2)+2\PriceLDEnergy[\SFn](t_1)
=& \intTslab{t_1}{t_2} 2t\chiLC(\dr\LDsubwt)(\dr \SFn)^2 -t\chiLC(\dr V_L)\LDsubwt|\dAng \SFn|^2\dFourPOneD \label{eLDTermsi}\\
&-\intTslab{t_1}{t_2} t\chiLC\frac{(\drrr\LDsubwt) \SFn^2}{2}\dFourPOneD \label{eLDTermsii}\\
&+\intTslab{t_1}{t_2} 2t\LDsubwt(\dr\chiLC)(\dr \SFn)^2\dFourPOneD \label{eLDTermsiii}\\
&-\intTslab{t_1}{t_2} t( 3(\dr\chiLC)(\drr\LDsubwt) +3(\drr\chiLC)(\dr\LDsubwt)  +(\drrr\chiLC)\LDsubwt)\frac{\SFn^2}{2}\dFourPOneD \label{eLDTermsiv}\\
&-\intTslab{t_1}{t_2} 2\dot{\SFn}\dot{g}(\dr \SFn)+\dot{\SFn}(\dr\dot{g})\SFn \dFourPOneD .\label{eLDTermsv}
\end{align}
The terms on the right in lines \eqref{eLDTermsiii}-\eqref{eLDTermsv} can be estimated by the local energy. (In these calculations, remember that $t^{-1}<\rs^{-1}$ and inverse powers of $t$ arise from differentiating $\chiLC$.)
\begin{align*}
\intTslice{t}& 2t\LDsubwt(\dr\chiLC)(\dr \SFn)^2 \dThreePOneD\\
<& C \intTLCslice{t} (\dr \SFn)^2\dThreePOneD < C \PriceMinEnergy ,\\
\intTslice{t}& |t( 3(\dr\chiLC)(\drr\LDsubwt) +3(\drr\chiLC)(\dr\LDsubwt)  +(\drrr\chiLC)\LDsubwt)\frac{\SFn^2}{2}| \dThreePOneD \\
<& C \intTLCslice{t} \frac{1}{1+|\rs|^{2}} |\SFn|^2 \dThreePOneD 
< C\PriceMinEnergy ,\\
\intTslice{t}& |2\dot{\SFn}\dot{g}(\dr \SFn)+\dot{\SFn}(\dr\dot{g})\SFn| \dThreePOneD\\
<& C \intTLCslice{t} |\dot{\SFn}|^2 +|\dr \SFn|^2 +\frac{1}{1+|\rs|^{2}} |\SFn|^2 \dThreePOneD\\
<&C \PriceMinEnergy . 
\end{align*}
The left-hand side can be estimated similarly by
\begin{align*}
\PriceLDEnergy(t)
<& tC \intTLCslice{t} |\dot{\SFn}|^2 +|\dr\SFn|^2 +\frac{1}{1+|\rs|^{2}} |\SFn|^2 \dThreePOneD \\
<&C t \PriceMinEnergy  .
\end{align*}
The two terms on the right appearing in line \eqref{eLDTermsi} are clearly positive, since $\LDsubwt$ is increasing and was chosen to go from negative to positive at the same value of $\rs$ as $-V_L'$. To control the term in \eqref{eLDTermsii} by the terms in \eqref{eLDTermsi}, we note that 
\begin{align*}
\drr\LDsubwt=& -\LDexp\LDpara\sgn(\rs)(1+\LDpara|\rs|)^{-(\LDexp+1)} ,\\
\drrr\LDsubwt
=& \LDexp(\LDexp+1)\LDpara^2 (1+\LDpara|\rs|)^{-(\LDexp+2)} -\LDexp\LDpara \delta(\rs) ,\\
\intTLCslice{t} -t\chiLC \frac{\drrr\LDsubwt}{2} |\SFn|^2 \dThreePOneD 
\geq& -t \intTLCslice{t} \chiLC \frac{\LDexp(\LDexp+1)\LDpara^2}{2} (1+\LDpara|\rs|)^{-(\LDexp+2)} |\SFn|^2 \dThreePOneD .
\end{align*}
We divide the range of integration into two pieces and use the estimate $t/|\rs| <C$ when $\rs>t/2$. From this, 
\begin{align*}
\intTslice{t} -t\chiLC \frac{\drrr\LDsubwt}{2} |\SFn|^2 \dThreePOneD 
\geq& -t \frac{\LDexp(\LDexp+1)\LDpara^2}{2} \intTLCtighterslice {t}(1+\LDpara|\rs|)^{-(\LDexp+2)} |\SFn|^2 \dThreePOneD \\
& - \frac{\LDexp(\LDexp+1)\LDpara^2}{2} \intTLCintermediateslice{t} (1+\LDpara|\rs|)^{-\LDexp+1} |\SFn|^2 \dThreePOneD .
\end{align*}
Applying the Hardy estimates, we find 
\begin{align*}
\intTslice{t} &-t\chiLC \frac{\drrr\LDsubwt}{2} |\SFn|^2 \dThreePOneD \\
\geq& -t \frac{C\LDexp(\LDexp+1)\LDpara^2}{2} \intTLCtighterslice{t} \chiLC\left( (\dr\LDsubwt)(\dr \SFn)^2 -(\dr V_L)\LDsubwt|\dAng \SFn|^2 \right) \dThreePOneD \\
&- C\PriceMinEnergy .
\end{align*}
Taking $\LDpara$ sufficiently small, we can dominate the integrand by half the terms in line \eqref{eLDTermsi}. Thus, 
\begin{align*}
\frac12 \intTslab{t_1}{t_2} t\chiLC\left(2(\dr\LDsubwt)(\dr \SFn)^2 -t(\dr V_L)\LDsubwt|\dAng \SFn|^2 \right)\dFourPOneD 
\leq& -2\PriceLDEnergy|_{t_1}^{t_2} + C \int_{t_1}^{t_2} \PriceMinEnergy(t) dt .
\end{align*}
The same estimate holds for $\Lie_{\RotElementi} \SFn$, so that summing over the components, we have
\begin{align*}
\frac12 \intTslab{t_1}{t_2} t\chiLC\left(2(\dr\LDsubwt)(\dr \dAng\SFn)^2 -t\chiLC(\dr V_L)\LDsubwt|\slap \SFn|^2\right) \dFourPOneD 
\leq& -2\PriceLDEnergy[\dAng\SFn]|_{t_1}^{t_2} + C \int_{t_1}^{t_2} \PriceMinEnergy[\dAng\SFn](t) dt .
\end{align*}
The left-hand side controls the trapping term by the Hardy estimate. Since the trapping term controls the growth of the conformal charge, 
\begin{align*}
\PriceConf[\SFn](t_2)
\leq& \PriceConf[\SFn](0)+ |-2\PriceLDEnergy[\dAng\SFn]|_{0}^{t_2}| + C \int_{0}^{t_2} \PriceMinEnergy[\dAng\SFn](t) dt ,\\
\leq& C \PriceConf[\SFn](0) +\sup_{t\in[0,t_2]}(t\PriceMinEnergy[\dAng\SFn](t)) + C \int_{0}^{t_2} \PriceMinEnergy[\dAng\SFn](t) dt .
\end{align*}
Since $\PriceEnergy[\dAng\SFn]$ is conserved, there is an immediate linear bound on the conformal charge. By applying the Cauchy-Schwartz estimate and integration by parts both twice, we can make the estimate
\begin{align*}
\PriceMinEnergy[\dAng\SFn](t) \leq& \PriceEnergy[\slap^2\SFn](t)^{1/4} \left(\frac{\PriceConf[\SFn](t)}{t^2}\right)^{3/4} .
\end{align*}
This allows us to make a self-improving estimate. From the linear bound, the conformal energy can't grow faster than $t^{1/4}$, and the $t^{1/4}$ implies a uniform bound. Thus, 
\begin{align*}
\PriceConf[\SFn](t)\leq& C (\PriceConf[\SFn](0) + \PriceEnergy[\slap^2\SFn](0)) , \\
\intTslabInfinity{0} t\chi |\dAng\SFn|^2 \dThreePOneD\leq& C (\PriceConf[\SFn](0) + \PriceEnergy[\slap^2\SFn](0)) .
\end{align*}

Applying the same argument with the factors of $t$ and $\chiLC$ dropped, so that $\LDwt(t,\rs)=\LDsubwt(\rs)$, we find
\begin{align*}
\intTslabInfinity{0} \frac{1}{(1+|\rs|^4)} |\SFn|^2 \dThreePOneD \leq& C \PriceEnergy[\SFn]. 
\end{align*}

%
%

\noindent{\bf \Large Acknowledgement}

This project was started when the author visited UC San Diego. The author would like to thank J. Sterbenz for his hospitality and for contributing several valuables ideas.

%
%

\bibliography{built}
\bibliographystyle{abbrv}

\end{document}